\documentclass[12pt]{article}
\usepackage{amsmath,amsfonts,amssymb}
\usepackage{graphicx}
\def\qed{\vbox{\hrule\hbox{\vrule height 1 ex\kern 1 ex\vrule}\hrule}}
\newcommand{\tmmathbf}[1]{\boldsymbol{#1}}
\newenvironment{proof}[1][Proof]{\begin{trivlist}{\item[]{\bf #1.\  %
}}}{{\small\hfill \qed} \end{trivlist}}

\newtheorem{theo}{Theorem}
\newtheorem{lem}[theo]{Lemma}
\newtheorem{prop}[theo]{Proposition}
\newtheorem{propri}[theo]{Properties}
\newtheorem{corol}[theo]{Corollary}
\newenvironment{defin}[1][Definition]{\begin{trivlist}{\item[]{\bf 
#1.\ }}}{ \end{trivlist}}

\def\thin{\preccurlyeq}
\def\Card{\mathrm{Card}}
\def\N{\mathbb{N}}
\def\Z{\mathbb{Z}}
\def\R{\mathbb{R}}

\def\rem{{\bf Remark. }}
\def\rems{{\bf Remark. }}
\def\cE{\mathcal{E}}

\def\cP{\mathcal{P}}
\def\l{\ell}
\def\eps{\varepsilon}

\def\un{\tmmathbf{1}}
\def\Ga{\mathrm{Gamma}}
				   
\begin{document}
\title{Products of Beta matrices and  sticky flows}
\author{Y. Le Jan, S. Lemaire}
\date{March 2004}
\maketitle
\begin{center}Universit\'e Paris-Sud\\Laboratoire de  Math\'ematique\\ B\^atiment
425\\ 91405 Orsay cedex\medskip\\ Yves.LeJan@math.u-psud.fr\\
Sophie.Lemaire@math.u-psud.fr
\end{center}
\vskip20pt
\begin{quote}{\small A discrete model of Brownian sticky flows on the unit circle is
described: it is constructed with products of Beta matrices on the discrete
torus.
 Sticky flows are defined by their ``moments'' which are
consistent systems of transition kernels on the unit circle.  Similarly, the moments of
the discrete model form  a  consistent system of transition matrices  on the discrete torus. A
convergence of Beta matrices to sticky kernels is shown at the level
of the moments. As the generators of the n-point processes are defined
in terms of Dirichlet forms, the proof is performed at the level of the
Dirichlet forms. The evolution of a probability measure by the flow of
Beta matrices is described by a measure-valued Markov process. A
convergence result of its finite dimensional distributions is deduced.
}
\end{quote}
\vskip20pt
\noindent{\bf Keywords~:} Markov chains
with continuous parameter, Polya urns, Dirichlet laws, stochastic flow of kernels,
Feller semigroups, Dirichlet forms, convergence of resolvents.
\medskip

\noindent{\bf 2000 MSC classifications~:} 60J27, 60J35, 60G09
\vskip20pt

\section{Introduction}
\label{intro}
In \cite{LeJanRaimond2,LeJanRaimond3}, a family of stochastic flows of
kernels on the circle 
$S^{1}$ called ``sticky flows'' is described. Sticky flows interpolate
between Arratia's flow of coalescing maps and the deterministic heat
flow. They are defined by their 
``moments'' which are consistent
systems of transition kernels on  $S^1$. In this article, a discrete 
version of sticky flows is presented for  sticky flows   
associated with  Brownian motions on $S^{1}$. This 
discrete model is defined by products of Beta matrices on the discrete 
torus $\frac{1}{N}(\Z/N\Z)$. It appears to be a
special case of a general construction which associates a flow of
Dirichlet matrices to any Markov chain on a finite set.\\ 
As in the  continuous case, 
the moments of the  flow of Beta matrices are consistent systems of transition matrices on 
$\frac{1}{N}(\Z/N\Z)$. The convergence of the flow of Beta matrices towards sticky
kernels is shown at the  level of the moments. Namely, it is
established that  for every $n\in
\N^*$, the semigroup
$\{P_{t}^{(n)},\ t\in\R\}$ of the $n$-point motion of the sticky flow
is approximated by the  $n$-point transition matrix of the 
flow of Beta matrices. Classical approximation theorems, such as theorem 6.1 in
\cite{EthierKurtz}, cannot be used since the generators of the
semigroups $P_{t}^{(n)}$ do not have a core of $C^2$ functions. The
generators are defined in terms of Dirichlet forms, thus the proof is
performed at the level of the Dirichlet forms. Given an initial law,
the flow of Beta matrices generates a measure-valued Markov process. A
convergence result of its finite dimensional distributions is
deduced.

Section 1 contains a  study of the flows of Dirichlet matrices
on a finite set. Flows of Beta matrices on $\frac{1}{N}(\Z/N\Z)$ and  Brownian sticky
flows on the circle $S^{1}$ are presented in section 2. Section 3 is
devoted to establishing the convergence of the flow of Beta matrices
on $\frac{1}{N}(\Z/N\Z)$
to a Brownian sticky flow on $S^{1}$.   
\section{Dirichlet matrices and Polya scheme}
A stochastic kernel on a finite set  $F$ is nothing but a random
transition matrix
on $F$. Given a sequence of i.i.d. random transition matrices  $(K_{i})_{i\in\Z}$ on
$F$, one may define a stochastic flow of kernels $(K_{s,t})_{ s\leq
  t}$ on $F$
 by setting for every $s\leq t$,
\begin{equation}
\label{defflot}K_{s,t}=\left\{
\begin{array}{ll}
 K_{s+1}\ldots K_{t-1}K_{t} &\mbox{ if the time is }\Z,\\
 K_{Z(s)+1}\ldots K_{Z(t)-1}K_{Z(t)}&\mbox{ if
  the time is } \R.\\
\end{array}\right.
\end{equation}where $Z$
  denotes a homogeneous Poisson process independent of $(K_i)_i$.\\
The sequence of matrices $(K_i)_i$ defines a random medium that evolves
with time. It is also a random distribution on paths. The $n$-point
motion is defined as an $n$-sample of this distribution. More
precisely, the environment (i.e. the sequence $(K_i)_{i}$) being fixed, one considers  the motion of $n$ independent inhomogeneous Markov
chains with transition kernels $(K_i)_i$. Under the annealed measure
(i.e. averaging on the environment),
the motion of $n$ points becomes  a homogeneous Markov chain with semigroup
${\{P_{t}^{(n)}:=E(K_{0,t}^{\otimes n})\}}$. The semigroups
$(P^{(n)}_t)_t$, $n\in\N^*$ form a consistent exchangeable family of Markovian
semigroups: all points play the same role and by
removing any point in the $n$-point motion, one gets the
$(n-1)$-point motion.    
In \cite{LeJanRaimond1}, following a generalization of  De
Finetti theorem, it  is shown that the law  of such a stochastic flow
of kernels is given by   a consistent
family of $n$-point Markovian semigroups 
$\{P_{t}^{(n)},\  n\in \N^*\}$.
\subsection{Dirichlet matrices}
We shall assume that the rows of every matrix $K_i$ are
independent Dirichlet random vectors.  A suitable choice of the
coefficients of the Dirichlet laws enables us to exhibit a consistent
and exchangeable system of probability measures $\{m^{(n)},\ n\in\N^*\}$ such that
for every $n\in \N^*$, $m^{(n)}$ is an invariant measure for the
$n$-point semigroup $(P_{t}^{(n)})_t$. \\
Let us first state a natural 
extension  of the Dirichlet distribution to the case of nonnegative
coefficients  and then define what
we shall call a Dirichlet random matrix and a flow of Dirichlet matrices. 
\begin{defin}
Let $\alpha_1,\ldots,\alpha_k$ be nonnegative reals such that at least one of them
are positive. Let  $1\leq i_1<\ldots<i_j\leq k$ denote the indices of the
positive coefficients and $V$ the set of points $x=(x_1,\ldots,x_k)\in[0,1]^{k}$ such that
$x_i=0$  if $i\not\in
\{i_1,\ldots,i_j\}$ and $\sum_{i=1}^{k}x_i=1$.\\
 A  random vector $X=(X_1,\ldots,X_k)$ is said to have the {\em
   Dirichlet distribution 
$\mathcal{D}(\alpha_1,\ldots,\alpha_k)$} if $P(X\not\in
V)=0$ and if $(X_{i_1},\ldots,X_{i_j})$ has the Dirichlet law 
$\mathcal{D}(\alpha_{i_1},\ldots,\alpha_{i_j})$ i.e. $$\frac{\Gamma(\alpha_{i_1}+\cdots+\alpha_{i_j})}{\Gamma(\alpha_{i_1})\cdots\Gamma(\alpha_{i_j})}x_{1}^{\alpha_{i_1}-1}\!\!\!\cdots
x_{j}^{\alpha_{i_j}-1}\!\!\tmmathbf{1}_{\{x_1+\cdots +x_j=1\}}dx_1\cdots dx_{j-1}.$$
\end{defin}
\begin{defin}Let $F$ and $G$ be two finite subsets. Let
  $A=(a_{i,j})_{(i,j)\in F\times G}$ be a matrix of
  nonnegative coefficients such that each row has at least a positive
  coefficient.\\
A random matrix   $X=(X_{i,j})_{(i,j)\in F\times G}$ will be called  {\em a
  Dirichlet matrix with parameter $A$} if the  rows of $X$ are
  independent random vectors and  the $i$-th row of $X$ has the Dirichlet
  distribution whose  parameters are given by the $i$-th row
  of $A$ for every $i\in F$. 
\end{defin}
\begin{defin}
Let $F$ be a finite set and $A=(a_{i,j})_{(i,j)\in F^2}$ be a matrix 
 of
  nonnegative coefficients such that each row has at least a positive
  coefficient.\\
The discrete time (or continuous time) stochastic flows defined by
 formula (\ref{defflot}) from  a family $(K_i)_{i\in \Z}$ of independent Dirichlet matrices of
 parameter $A$   is said to
  be  {\em a  stochastic flow of Dirichlet matrices of parameter $A$
 on $F$}.  
\end{defin}
Let us notice that  the product of two independent Dirichlet
matrices is a matrix whose rows are generally neither independent nor distributed according
to a Dirichlet distribution. Still,  
the following lemma holds; it states a probably well-known property of
Dirichlet laws which  is
based on properties of Gamma distributions (first?)  mentioned by
Edwin Pitman in \cite{EPitman}. 
\begin{lem}
\label{lemgamma}
Let $r$ and $s$ be positive integers. Let  $A=(a_{i,j})_{1\leq i\leq
 r,\ 1\leq j\leq s}$ be a matrix 
 of
  nonnegative coefficients such that each row has at least a positive
  coefficient.
Let  $X$ be a Dirichlet  matrix of parameter $A$.\\
  If
$Y$ is  a random vector independent of the random matrix $X$,
 with Dirichlet distribution
$\mathcal{D}(\sum_{j=1}^{s}a_{1,j},\ldots,\sum_{j=1}^{s}a_{r,j})$, 
then the random vector $YX$  has the Dirichlet distribution
$\mathcal{D}(\sum_{i=1}^{r}a_{i,1},\ldots,\sum_{i=1}^{r}a_{i,s})$.
\end{lem}
\begin{proof}
We shall construct copies of $X$ and $Y$ using a family of independent
Gamma random variables. This construction is based on the following
properties of the Gamma laws:
\begin{propri}
\label{properGam}Let $\alpha_1,\ldots,\alpha_r$ be nonnegative reals and
  $Z_1,\ldots,Z_r$ be a family of independent random variables with
  Gamma laws   $\Ga(\alpha_1,1),\ldots,$ $\Ga(\alpha_r,1)$ respectively (by
  convention, the Gamma law of parameter $0$ is the Dirac at
  $0$). Assume that some of the parameters $\alpha_1,\ldots,\alpha_r$
  are positive. Then
\begin{itemize}
\item[-]the random variable $Z=\sum_{i=1}^{r}Z_i$  follows
  the law $\Ga(\sum_{i=1}^{r}\alpha_i,1)$.
\item[-] the random vector
  $(\frac{Z_1}{Z},\ldots,\frac{Z_r}{Z})$
  has the Dirichlet law $\mathcal{D}(\alpha_1,\ldots,\alpha_r)$ and is
  independent of the random variable $Z$. 
\end{itemize}
\end {propri} 
These properties can be shown by computing the Laplace transform of
the random vector $(\frac{Z_1}{Z}, \ldots,\frac{Z_r}{Z},Z)$. \\
Let $\{U_{i,j},\ (i,j)\in\{1,\ldots,r\}\times\{1,\ldots,s\}\}$ denote a family
of independent random variables such that $U_{i,j}=0$ if $a_{i,j}=0$
and  $U_{i,j}$ has the
$\Ga(a_{i,j},1)$ distribution if $a_{i,j}>0$. Set $U_{i}=\sum_{j=1}^{s}U_{i,j}$, $U=\sum_{i=1}^{r}U_i$ and
$V_{i,j}=\frac{U_{i,j}}{U_{i}}$ for every $(i,j)\in
\{1,\ldots,r\}\times\{1,\ldots,s\}$. It follows from the recalled properties of the Gamma laws that for
every $i\in \{1,\ldots,r\}$,  
  the random vector $V_i\!=\!(V_{i,1},\ldots,V_{i,s})$ has the Dirichlet
  distribution $\mathcal{D}(a_{i,1},\ldots,a_{i,s})$  and  is independent of $U_i$ that
  has the $\Ga(\sum_{j=1}^{s}a_{i,j},1)$ distribution. 
  As ${\{U_{i,j},\ (i,j)\in\{1,\ldots,r\}\times\{1,\ldots,s\}\}}$ is a family of
independent random variables, the random variables
$V_1,\ldots,V_r,U_1,\ldots,U_r$ are independent. In particular,
  $V=(V_{i,j})_{i,j}$ is a Dirichlet matrix of parameter $A$.\\
It also follows that  the random vector
$W=(\frac{U_{1}}{U},\ldots,\frac{U_{r}}{U})$ has the same law as $Y$,
i.e. the Dirichlet distribution 
$\mathcal{D}(\sum_{j=1}^{s}a_{1,j},\ldots,\sum_{j=1}^{s}a_{r,j})$  and is independent
of the random vectors $V_{1},\ldots,V_{r}$.\\
If we set $T_j=\sum_{i=1}^{r}W_iV_{i,j}$ for $j \in\{1,\ldots,s\}$, then
$T_j=\frac{1}{U}\sum_{i=1}^{r}U_{i,j}$. As the random variables
$\sum_{i=1}^{r}U_{i,j}$ for $j\in\{1,\ldots,s\}$ are
independent and have the
$\Ga(\sum_{i=1}^{r}a_{i,1},1)$,\ldots, $\Ga(\sum_{i=1}^{r}a_{i,s},1)$ distributions
respectively, the random vector $(T_1,\ldots,T_s)$ has the Dirichlet
distribution $\mathcal{D}(\sum_{i=1}^{r}a_{i,1},\ldots,\sum_{i=1}^{r}a_{i,s})$.
 \qed\end{proof}
\subsection{Invariant distributions}
Let us notice that if $P$ is a stochastic matrix indexed by $F$ and $m$ is a
finite positive measure on $F$ invariant by $P$ then the matrix $A=(a_{i,j})$
defined  by $a_{i,j}=m_{i}P_{i,j}$ for every $(i,j)\in F^2$ is such
that the  sum of the coefficients of the $i$-th row  is equal to the sum of
the coefficients of the $i$-th
column  for every $i\in F$. Thus applying  lemma
\ref{lemgamma} yields:
\begin{prop}\label{mesinvnoy}   
Let $P$ be a stochastic matrix indexed by a finite set $F$ and let  $m$ be
a positive measure on $F$  invariant by $P$. Set $A=(a_{i,j})$ the matrix 
defined by  $a_{i,j}=m_ip_{i,j}$ for every $(i,j)\in F^{2}$. Let 
$(K_{s,t})_{s\leq t}$ be  a discrete time stochastic flow of Dirichlet 
matrices of parameter $A$. Let $\mu$ be a Dirichlet vector on $F$ of 
parameter $m$, independent of  $(K_{s,t})_{s\leq t}$.
\begin{enumerate}
    \item[(i)] For every $ s\leq t$, $\mu K_{s,t}$ is 
a Dirichlet vector of parameter $m$.
\item[(ii)] For every $n\in \N^{*}$, $E(\mu^{\otimes n})$ is an 
invariant probability measure for the $n$-point semigroup 
$(P_{t}^{(n)})_{t}$ associated with the Dirichlet flow $(K_{s,t})_{s\leq t}$. Its expression can be given iteratively as
follows:
\begin{equation}\label{recmesinv}
E(\mu^{\otimes n})(x)=\frac{m(x_n)+\sum_{i=1}^{n-1}\un_{\{x_n=x_i\}}}{m(F)+n-1}
E(\mu^{\otimes (n-1)})(\underline{x}).
\end{equation}
for every $x=(\underline{x},x_{n})\in F^n$.
\item[(iii)] Assume that $P$ is periodic with period $d>1$. There 
exists a partition $C_{0},\ldots,C_{d-1}$ of $F$ such that $p_{i,j}>0$ 
only if there exists $r\in\{0,\ldots,d-1\}$ such that $i\in C_{r}$ and 
$j\in C_{r+1}$ (by convention $C_{jd+k}=C_{k}$ for every $j\in \N$ 
and $k\in \{0,\ldots,d-1\}$). Then for every 
$k,r\in\{0,\ldots,d-1\}$ and $j\in \N$, $\mu_{\mid C_{k}}K_{0,dj+r}$ is a 
Dirichlet vector on $F$ of parameter $m\un_{C_{k+r}}$. 
\end{enumerate}
\end{prop}
\begin{proof}
Assertion $(i)$ follows from  lemma \ref{lemgamma}. For every
    $n\in \N^*$ and $t\in \N$, $E((\mu K_{0,t})^{\otimes
    n})$ is equal to $E(\mu^{\otimes n})P_{t}^{(n)}$, thus $E(\mu^{\otimes n})$ is
   invariant by $P_{t}^{(n)}$. \\
Let us compute the moments of $\mu$.  Set $e_{1},\ldots,e_r$ the elements of
    $F$. For $x\in F^{n}$, let $s_{e}(x)$ denote the number of
    coordinates of $x$ that equal to $e$. As  $E(\mu^{\otimes
    n})(x)$ is the moment of order $(s_{e_{1}}(x),\ldots,s_{e_{r}}(x))$ of
    the Dirichlet law $\mathcal{D}(m(e_1),\ldots,m(e_r))$, 
\begin{eqnarray*}
E(\mu^{\otimes
    n})(x)&=&\frac{\Gamma(m(F))}{\Gamma(m(F)+n)}
\prod_{e\in
    F}\frac{\Gamma(m(e)+s_{e}(x))}{\Gamma(m(e))}\\
&=&\frac{1}{\gamma_{m(F)}(n)}\prod_{e\in
    F}\gamma_{m(e)}(s_{e}(x)).
\end{eqnarray*}
where $\gamma_{a}(u)$ denotes the product $\prod_{i=0}^{u-1}(a+i)$ for
    every $u\in \N^*$ and $a\in \R$. To obtain the iterative
    expression, it remains to notice that if $x=(\underline{x},x_n)$ then  
    $s_{e}(x)=\un_{\{x_n=e\}}+s_{e}(\underline{x})$. \\ 
 The
    restriction of $\mu$ to $C_k$ is the random measure
    $\mu_{|C_k}(\ \cdot\ )=\frac{\mu(C_k\cap\ \cdot\ )}{\mu(C_k)}$. As $\mu$ is a
    Dirichlet vector of parameter $m$,  $\mu_{|C_k}$ is a Dirichlet
    vector of parameter $m\un_{C_k}$. Lemma
    \ref{lemgamma} applies to the vector $(\mu_{|C_{k}}(i))_{i\in
    C_k}$ and the submatrix $\tilde{K}_l=(K_{l,l+1}(i,j))_{i\in C_k,\ j\in
    C_{k+1}}$ of $K_{l,l+1}$ for every $l\in\Z$;  this shows that
    $\mu_{|C_{k}}K_{l,l+1}$ is a Dirichlet vector of parameter
    $m\un_{C_{k+1}}$ for every $k\in \{0,\ldots,d\}$ and $l\in \Z$.
 \qed\end{proof}
In order to present an explicit expression of $E(\mu^{\otimes n})$, let us
introduce some notations associated with a partition 
 $\pi$ of $[n]=\{1,\ldots,n\}$.\\
Let $\mathcal{P}_n$ be the set of
all partitions of $[n]$. Let $|\pi|$ 
denote the number of  non-empty blocks of $\pi$. The symbol $\thin$
will be used to designate an order relation between partitions:
$\pi'\thin \pi$ if $\pi'$ is equal to $\pi$ or is a finer partition than $\pi$. Let 
$C_{\pi}$ 
be the set of points $x\in F^n$ such that 
$x_{i}=x_{j}$ if and only if $i$ and $j$ are in the same block of 
$\pi$. Set $E_{\pi}=\cup_{\pi',\ \pi\thin \pi'}C_{\pi'}$. It is  
the set of points $x\in F^n$ such that, 
if $i$ and $j$ are in the same block of 
$\pi$, then $x_{i}=x_{j}$. It  is isomorphic to
$F^{|\pi|}$. Let $\phi_{\pi}:F^{|\pi|}\mapsto E_{\pi}$ be the  one-to-one
mapping  defined   as
follows:
for all $x\in F^{|\pi|}$, $\phi_{\pi}(x)=(y_{1},\ldots,y_{n})$  where 
$y_{i}=x_{l}$ if $i$ belongs to the $l$-th block of $\pi$. 
The measure $E(\mu^{\otimes n})$ can be expressed as a 
combination of  probability measures on the sets $E_{\pi}$, $\pi \in \mathcal{P}_n$ as follows: 
\begin{prop}
\label{decommesurinv}Let $m$ be a measure on a finite  set $F$ and
let $\mu$ be a Dirichlet vector on   $F$ with parameter
$m$. Set $\tilde{m}=\frac{m}{m(F)}$. For every $n\in \N^*$, 
\begin{equation}
\label{exprecmesinv}
E(\mu^{\otimes n})=\sum_{\pi\in
  \mathcal{P}_n}p^{(m(F))}_{\pi}\phi_{\pi}(\tilde{m}^{\otimes |\pi|})
\end{equation}
where
$p^{(a)}_{\pi}=a^{k}\frac{\overset{k}{\underset{i=1}{\prod}}(n_i-1)!}{\overset{n-1}{\underset{i=0}{\prod}}(a+i)}$
if $\pi$ is a partition of $[n]$ with $k$ nonempty blocks of length
$n_1,\dots,n_k$ and $a$ is a positive real.  
\end{prop}
\begin{proof} The proof can be deduced from the  expression
  (\ref{recmesinv}) of $E(\mu^{\otimes n})$,  by induction  on $n$. Assume
   that formula (\ref{exprecmesinv}) is true for
  $E(\mu^{\otimes (n-1)})$. Fix a
  point  $x=(\underline{x},x_n)\in F^n$. If $\pi$ is the partition of
   $[n-1]$ such that $\underline{x}\in C_{\pi}$, then 
$$E(\mu^{\otimes n})(x)=\sum_{\pi',\ \pi'\thin
  \pi}p_{\pi'}^{(m(F))}\phi_{\pi'}(\tilde{m}^{\otimes
  |\pi'|})(\underline{x})\frac{m(x_n)+\sum_{i=1}^{n-1}\un_{\{x_i=x_n\}}}{m(F)+n-1}.$$ 
Set $B_1,\ldots,B_k$ denote the nonempty blocks of the partition
  $\pi$. Let  $s(\pi)$ denote the partition $(\pi,\{n\})$ and $s_{i}(\pi)$
  denote the partition of $[n]$ obtained from $\pi$ by adding $n$ to
  the block $B_i$  for $i\in\{1,\ldots,k\}$. Two cases arise
  according as $x\in C_{s(\pi)}$ or $x\in C_{s_{i}(\pi)}$ for some
  $i\in\{1,\ldots,k\}$ since for every $a>0$, $$p^{(a)}_{s(\pi)}=
\frac{a}{a+n-1}p^{(a)}_{\pi} \mbox{ and }p^{(a)}_{s_{i}(\pi)}=
\frac{\Card(B_i)}{a+n-1}p^{(a)}_{\pi} \mbox{ for }i\in\{1,\ldots,k\}.$$ Then  
   formula (\ref{exprecmesinv}) for $E(\mu^{\otimes n})$ follows from
  the fact that:
\begin{itemize}
\item[-]
   the partitions $\tilde{\pi}$ finer than $s(\pi)$ have the form
$s(\pi')$ where $\pi'$ is finer than $\pi$,
\item[-] the partitions finer than $s_i(\pi)$ either have the
  form $s(\pi')$ where $\pi'$ is a partition of $[n-1]$ finer than
  $\pi$,  or have the form $s_j(\pi')$ where  $\pi'$ is a partition of $[n-1]$ finer than
  $\pi$ such that its $j$-th block is a subset of $B_i$.
\end{itemize}
 \end{proof}
\rem $(p^{(a)}_{\pi})_{\pi\in\mathcal{P}_n}$ is the exchangeable
partition function of an exchangeable sequence of random variables
governed by the
Blackwell-MacQueen Urn scheme \cite{Pitmanurn}. 
\subsection{Description of the $n$-point motion}
In order to complete the description of the $n$-point motion
associated with the discrete time  flow of Dirichlet matrices of parameter $A$
denoted by $(K_{s,t})_{s\leq t}$, let us
compute the transition matrix $P^{(n)}$: for
$x=(x_1,\ldots,x_n)$ and $y=(y_1,\ldots,y_n)$ in
$F^n$, $$P^{(n)}(x,y)=E(\prod_{i=1}^{n}K_{1}(x_i,y_i))=
\prod_{l\in F}E(\prod_{h\in F}K_{1}(l,h)^{s_{l,h}(x,y)})$$ where
$s_{l,h}(x,y)=\Card(\{i,\ (x_i,y_i)=(l,h)\})$. Let us set
$\gamma_a(u)=\prod_{i=0}^{u-1}(a+i)$ and
$s_l(x)=\Card(\{i\in\{1,\ldots,n\},\ x_i=l\})$. It follows from the
expression of the moments of the Dirichlet laws that
$$P^{(n)}(x,y)=\prod_{i\in F}\frac{\prod_{j\in F}
\gamma_{a_{i,j}}(s_{i,j}(x,y))}{\gamma_{\sum_{j\in F}a_{i,j}}(s_{i}(x))}.$$ Let us
remark that:
$$
 P^{(1)}(x,y)=\frac{a_{x,y}}{\sum_{z\in F}a_{x,z}}   \mbox{ and }
 P^{(n)}(x,y)=\frac{a_{x_n,y_n}+s_{x_n,y_n}(\underline{x},\underline{y})}
 {\sum_{j\in
 F}a_{x_n,j}+s_{x_n}(\underline{x})}P^{(n-1)}(\underline{x},\underline{y})$$
 for $n\geq 2$ 
 where $\underline{x}=(x_1,\ldots,x_{n-1})$.\\
 Thus the transition mechanism can be described as follows: the first
point moves from a site $i$ to a site $j$ with probability
$p_{i,j}=\frac{a_{i,j}}{\sum_{\l\in F}a_{i,\l}}$.  The motion of the $(k -1)$ first points being known,
the $k$-th point moves from a site $i$ to a site $j$ with
probability $\frac{a_{i,j}+u}{\sum_{\l\in F}a_{i,\l}+ v}$ if among the $k - 1$ first
points, $v$ were located on the site $i$, and $u$ have moved from
$i$ to the site $j$. This is a combination of independent Polya
urns attached at each site.\\

\noindent Let us deduce some elementary properties of $P^{(n)}$:
\begin{enumerate}
\item
If the one-point motion is defined as a reversible Markov chain on $F$, 
 then for every $n\in \N^*$ the $n$-point motion  is a reversible
 Markov chain on $F^n$.
\item If $A$ is irreducible and aperiodic then for every $n\in \N^*$,
  $P^{(n)}$ is also an irreducible and aperiodic matrix.
\item Assume that $A$ is irreducible and periodic of period $d$. Let   
  $C_{0},\ldots,C_{d-1}$ denote a partition of $F$ such that for every
  $r\in\{0,\ldots,d-1\}$,  $i\in C_r$ and $j\not\in C_{r+1}$ imply
  $a_{i,j}=0$. Then for every $n\in\N^*$ and $r\in\{0,\ldots,d-1\}$, 
$(C_{r})^{n}$ is a closed subset for $P^{(n)}$ and 
  $(P^{(n)})^d$ is an irreducible aperiodic matrix on $(C_{r})^{n}$.
\end{enumerate}
The property (1) implies that if the one-point motion is a reversible
Markov chain on $F$ with reversible probability measure $m$ and if $\mu$ is a
Dirichlet vector on $F$ with parameter $m$ independent of the stochastic
kernels $(K_i)_i$, then the stationary vector-valued Markov process
$(\mu K_{0,j})_j$ is reversible.  
The properties (2) and (3) of the $n$-point motions and  proposition
 \ref{mesinvnoy} enable us to establish the asymptotic behaviour of the
  flow of Dirichlet matrices:  
\begin{corol}
\label{asympmatrices}   
Let $P$ be a stochastic matrix indexed by a finite set $F$ and let  $m$ be
a positive measure on $F$  invariant by $P$. Set $A=(a_{i,j})$ the matrix 
defined by  $a_{i,j}=m_ip_{i,j}$ for every $(i,j)\in F^{2}$. Let 
$(K_{s,t})_{s\leq t}$ be  a discrete time stochastic flow of Dirichlet 
matrices of parameter $A$. 
\begin{itemize}
\item[(i)]Assume that $P$ is irreducible and aperiodic.\\ For every
  probability measure $\nu$ on $F$, $\nu K_{0,j}$ converges in
  law to a Dirichlet vector of parameter $m$ as $j$ tends to
  $+\infty$. 
\item[(ii)]Assume that $P$ is irreducible and periodic of period
  $d$. Let $C_{0},\ldots,C_{d-1}$ denote a partition of $F$ such that for every
  $r\in\{0,\ldots,d-1\}$, if $i\in C_{r}$  and  $j\not\in C_{r+1}$
  then $p_{i,j}=0$. Set $C_{k+d}=C_{k}$ for every $k\in \N$. For $k\in\{0,\ldots,d-1\}$, let $\nu_{k}$ denote a
  probability measure on $F$ such that $\nu_k(C_k)=1$. \\
For every $r$ in $\{0,\ldots,d-1\}$, 
  $(\nu_{0}K_{0,jd+r},\ldots,\nu_{d-1}K_{0,jd+r})$ converges in law,
  as $j$ tends to $+\infty$,  to a 
  vector $(\mu_r,\ldots,\mu_{r+d-1})$ of independent Dirichlet
  vectors  of parameters
  $m\un_{C_{r}},\ldots,m\un_{C_{r+d-1}}$ respectively. 
 \end{itemize}
\end{corol}
\begin{proof}
\begin{itemize}
\item[(i)] $(\nu K_{0,j})_j$ is a sequence of random variables that
  take their values in the set of probability measures on $F$. To
  prove that they converge in law to a random measure $\mu$, it
  suffices to prove that for every $n\in\N^*$, the sequence of
  probability measures $(E((\nu K_{0,j})^{\otimes n}))_j$ on $F^{n}$
  converges  to $E(\mu^{\otimes n})$. As
  $P^{(n)}_{0,j}=E(K_{0,j})^{\otimes n}$ is the transition matrix of
  an irreducible aperiodic Markov chain on $F^{n}$, $(\nu^{\otimes
  n}P_{0,j}^{(n)})_j$ converges  to the stationary law of the
  Markov chain that is $E(\mu^{\otimes n})$ where $\mu$ is a Dirichlet
  vector of parameter $m$, by proposition \ref{mesinvnoy}.  
\item[(ii)]Let us first note that for every $j\in \N$,
  $\nu_0K_{0,j},\ldots,\nu_{d-1}K_{0,j}$ are independent variables
  since $K_{0,j}$ is a product of $j$ independent random matrices, each
  of them having independent rows. As $(P^{(n)})^{d}$ is an irreducible and aperiodic
  Markov chain on $(C_{u})^{n}$ for every $u\in\{0,\ldots,d-1\}$ and
  $n\in \N^{*}$, it
  can be shown as in the proof of $(i)$, that for every
  $k,r\in\{0,\ldots,d-1\}$, $\nu_kK_{0,dj+r}$ converges in law to a
  Dirichlet vector of parameter $m\un_{C_{k+r}}$ as $j$ tends to
  $+\infty$.
\end{itemize}
 \end{proof}
\section{A discrete model of a sticky flow on $S^{1}$\label{sectbetaflow}}
In the remaining part of this paper, we shall consider  a particular   flow
of Dirichlet matrices on the lattice $T_N=\frac{1}{2N}(\Z/2N\Z)$  of $S^{1}$,
called {\em flow of Beta matrices}. It  gives  a simple discrete model
of the Brownian sticky
flow defined by Y. Le Jan and
O. Raimond in \cite{LeJanRaimond2,LeJanRaimond3}. After a description of  the flow
of Beta matrices, the main properties of the Brownian sticky flows
will be summarized. 
\subsection{Beta matrices}
 Let  $a$ be 
a positive real, $N\in \N^*$  and $K$ be  a  matrix  on the
lattice $T_N$ defined as follows:
  $$K ( i, j ) = X_i \tmmathbf{1}_{\{j = i + \frac{1}{2N}\}} + ( 1
- X_i ) \tmmathbf{1}_{\{j = i - \frac{1}{2N}\}}\mbox{ for every } i,j\in T_N$$ where $X_1, \ldots, X_{2N}$
are independent Beta$( \frac{a}{2N}, \frac{a}{2N} )$ random
variables.\\
The matrix $K$ is a Dirichlet matrix of parameter $A=(m_ip_{i,j})_{i,j}$ where
$P=(p_{i,j})_{i,j}$ is the transition matrix of the symmetric random
walk on $T_N$ and $m$ is the uniform measure on $T_N$ with total mass
$2a$.  
 Let $(K_{i})_{i\in \Z}$ be a sequence
of independent matrices having the same law as $K$ and let $Z$  be an independent Poisson process on $\mathbb{R}$
with intensity $4N^2$. The continuous time Dirichlet
flow of matrices defined on $T_N$ by
$$K_{N,s,t}=K_{Z(s)+1}K_{Z(s)+2}\cdots K_{Z(t)-1}K_{Z(t)} \mbox{ for every }s\leq t,$$
will be called  {\em a  flow of Beta matrices on $T_N$ of parameter
  $a$}. Let us note that $A$ is
irreducible and two-periodic thus the $n$-point process is not
irreducible. We shall focus our attention on the following
irreducible set: 
\begin{eqnarray*}
T_{N}^{(n)}&=&\{\frac{x}{2N},\ x\in(\Z/2N\Z)^n\mbox{ and }x_i-x_1\in 2\Z , \mbox{ for all }1\leq
i\leq n\}\\
&=&(\frac{1}{N}(\Z/N\Z))^{n}\cup(\frac{1}{2N}((2\Z+1)/2N\Z))^{n}
\end{eqnarray*}
The  transition matrix of the $n$-point motion  on $T_{N}^{(n)}$ at time $t$ will be
denoted $P^{(n)}_{N,t}$. The jump Markov chain of the $n$-point
process is two-periodic: in one step, a point with odd coordinates moves to a
point with even coordinates and conversely. Its  transition matrix
will be denoted $P^{(n)}_N$. It has a reversible  probability measure
denoted by
$m^{(n)}_N$:  
\begin{equation}
\label{mesinvbeta}
m^{(n)}_N=\frac{1}{2}(E(\mu_{N}^{(0)\otimes
n})+E(\mu_{N}^{(1)\otimes n}))
\end{equation}
where $\mu_{N}^{(0)}$ and $\mu_{N}^{(1)}$ denote
  independent 
  Dirichlet vectors on $T_N$ with parameter
  $m\un_{\frac{1}{N}\Z/N\Z}$ and $m\un_{\frac{1}{2N}((2\Z +1)/2N\Z)}$
    respectively. Let $\eta$ be a Bernoulli random variable with
  parameter $\frac{1}{2}$, independent of $\mu_{N}^{(0)}$
    and $\mu_{N}^{(1)}$. It follows that if $\nu$ is a probability measure on $\frac{1}{N}\Z/N\Z$
or $\frac{1}{2N}((2\Z+1)/2N\Z)$ then $(\nu K_{N,0,t})_{t\in\R_+}$ converges
in law to $\mu_{N}^{(\eta)}$ as $t$ tends to $+\infty$.\\

\begin{figure}[!p]
\parbox{2.9in}{%
\includegraphics*[scale=0.4]{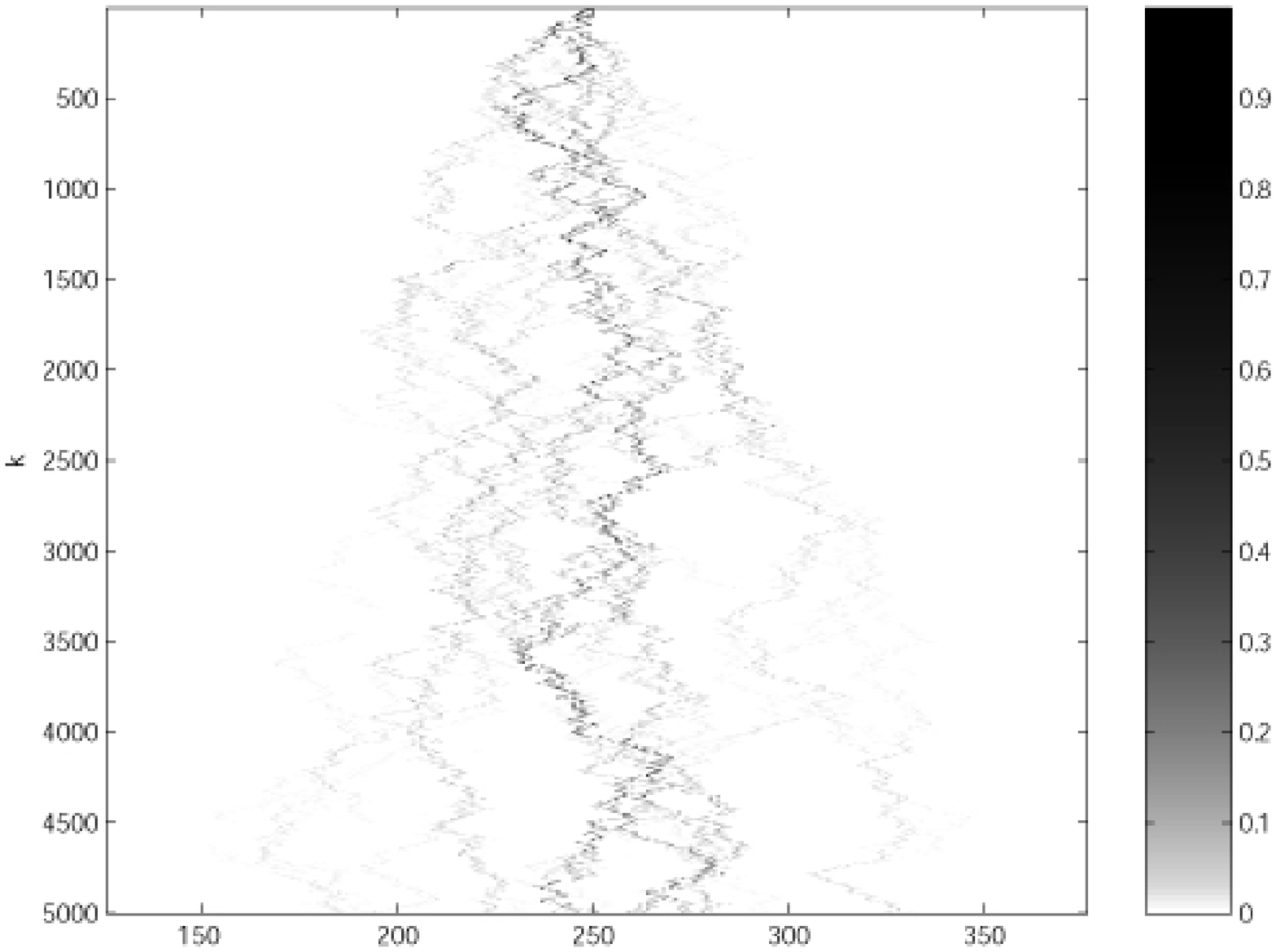}
\caption{a sample of the measure-valued process $\{\nu^{(20)}_k,\ 0\leq k\leq  5000\}$}
\label{fignoyau20}}\qquad
\begin{minipage}{2.9in}
\includegraphics*[scale=0.4]{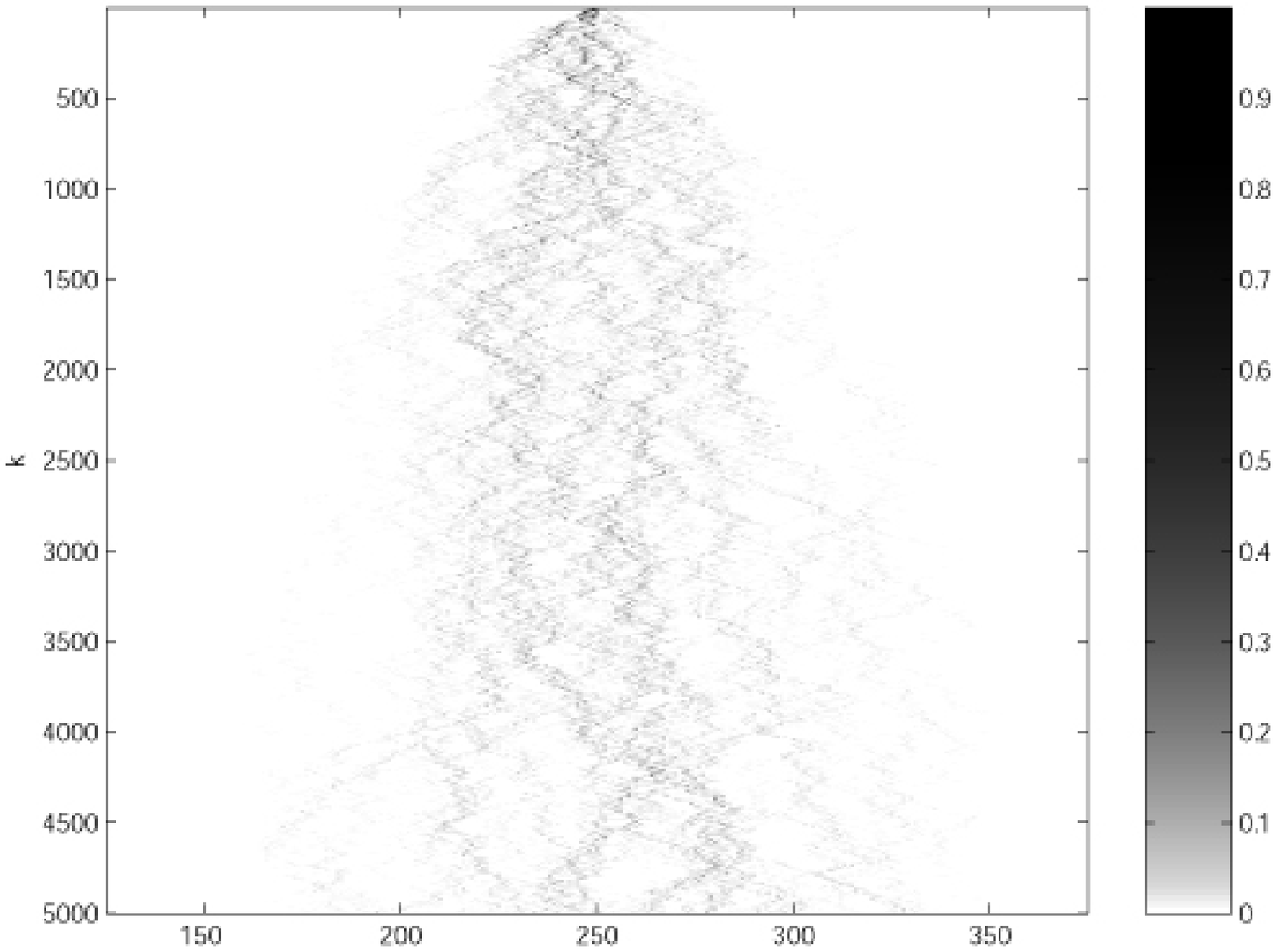}
\caption{a sample of the measure-valued process $\{\nu^{(100)}_{k},\ 0\leq k\leq  5000\}$}
\label{fignoyau100}
\end{minipage}
\end{figure}

\begin{figure}[!p]\parbox{2.9in}{%
\includegraphics*[scale=0.4]{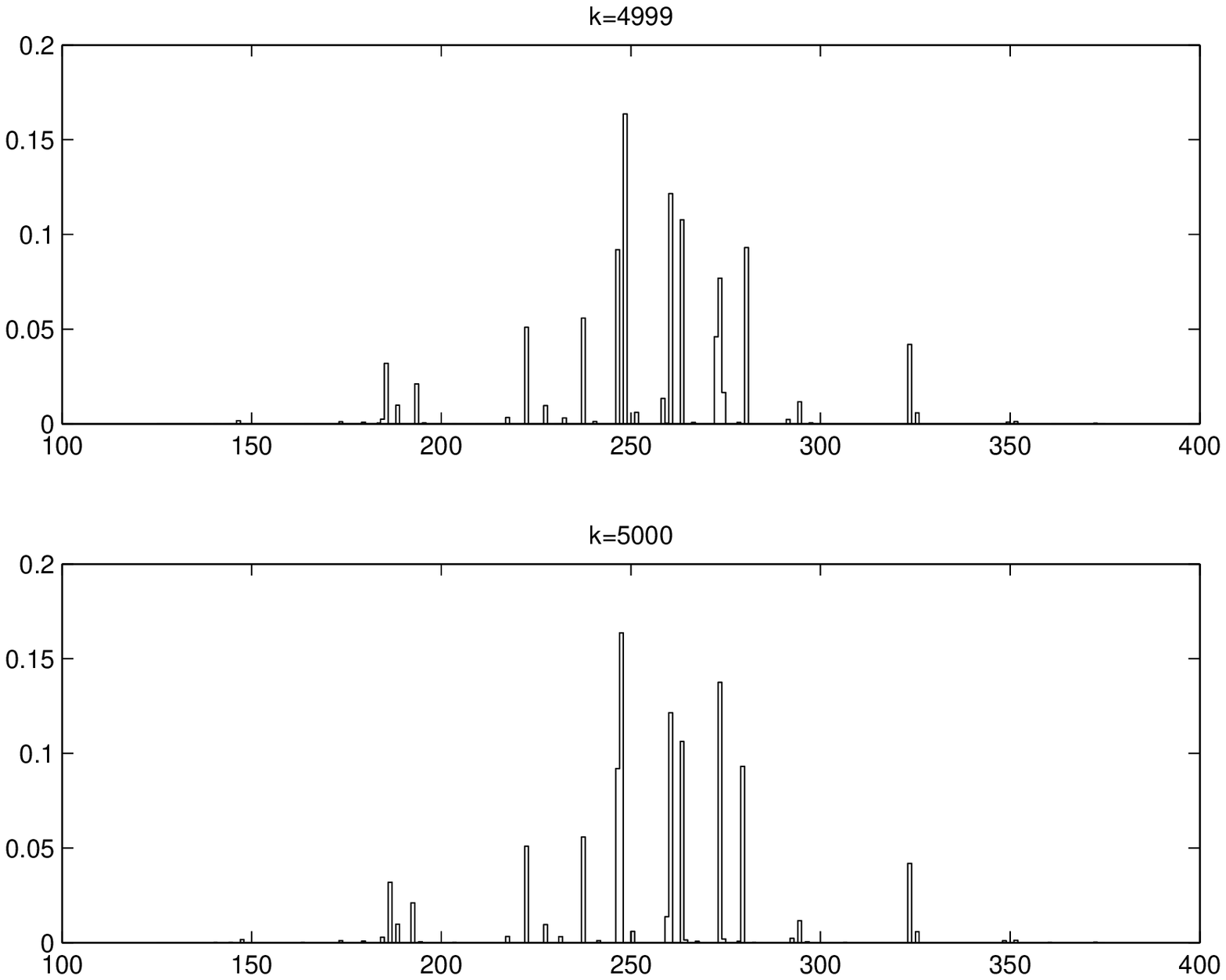}
\caption{histograms of the sample of the probability measure
  $\nu^{(20)}_{k}$ for  two successive values of $k$: $k=4999$ and $k=5000$}
\label{histo20}}\qquad
\begin{minipage}{2.9in}
\includegraphics*[scale=0.4]{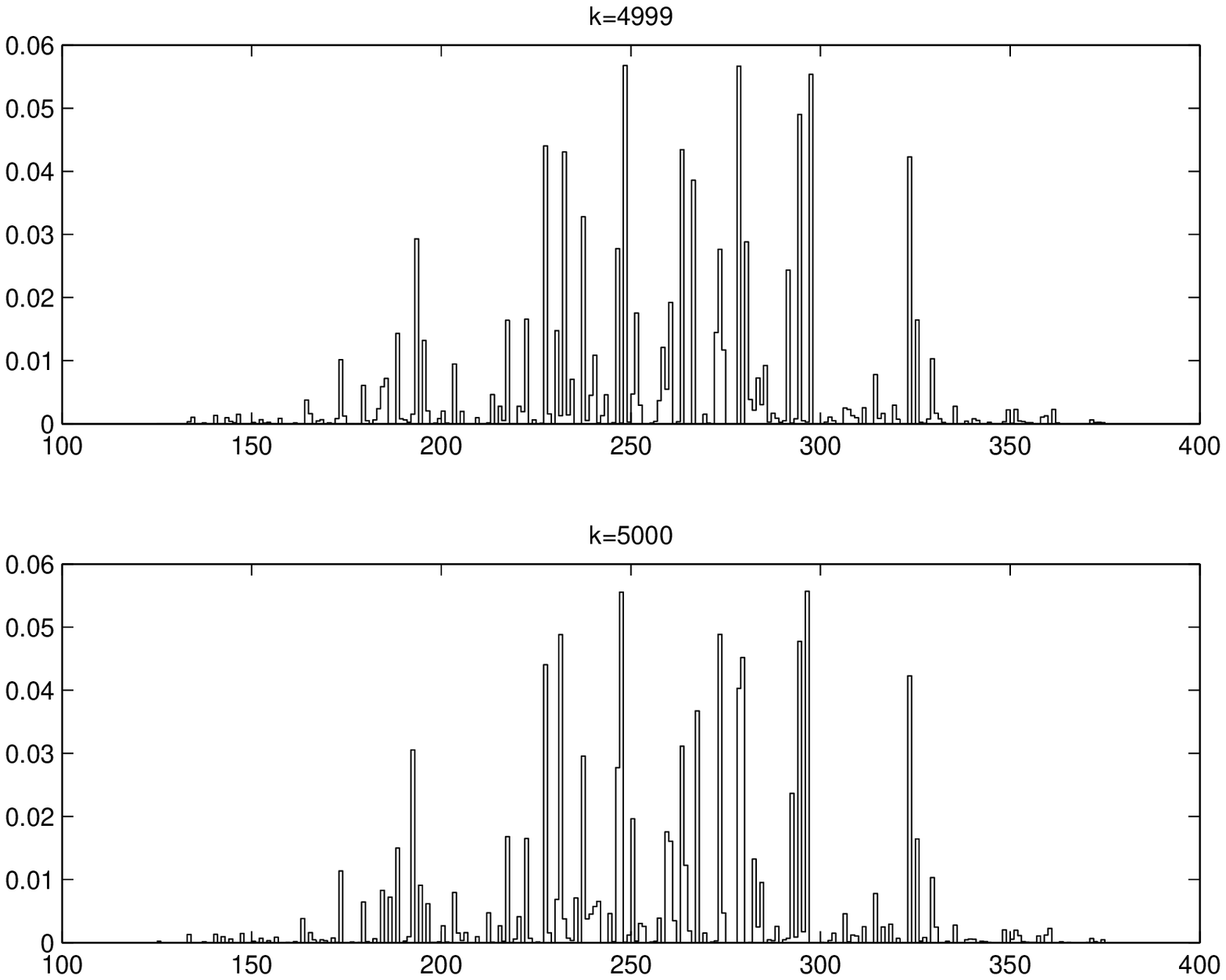}
\caption{histograms of the sample of the probability measure
  $\nu^{(100)}_{k}$ for  two successive values of $k$: $k=4999$ and $k=5000$}
\label{histo100}
\end{minipage}
\end{figure}

\noindent A sample path of the measured-valued Markov process 
  $\nu^{(a)}\!=\!(\delta_{x}K_{N,0,2k})_{k\in\N}$  associated to the discrete-time
  flow of Beta matrices on
  $\frac{1}{N}(\Z/N\Z)$, is represented in figures \ref{fignoyau20} and
  \ref{fignoyau100} for  $N=500$, $x=\frac{1}{2}$ and two choices of the parameter $a$:
  $a=20$ and $a=100$.
  In these figures, the sample of
  $\nu^{(a)}_k$ at a fixed $k$ is represented by a horizontal
  line of $N$ colored  pixels. A detail of the histogram of the sample of
  $\nu^{(a)}_k$ for $k=500$ and $k=1000$ is given figures
  \ref{histo20}
and \ref{histo100}.
\begin{figure}[!p]
\parbox{2.9in}{%
\includegraphics*[scale=0.4]{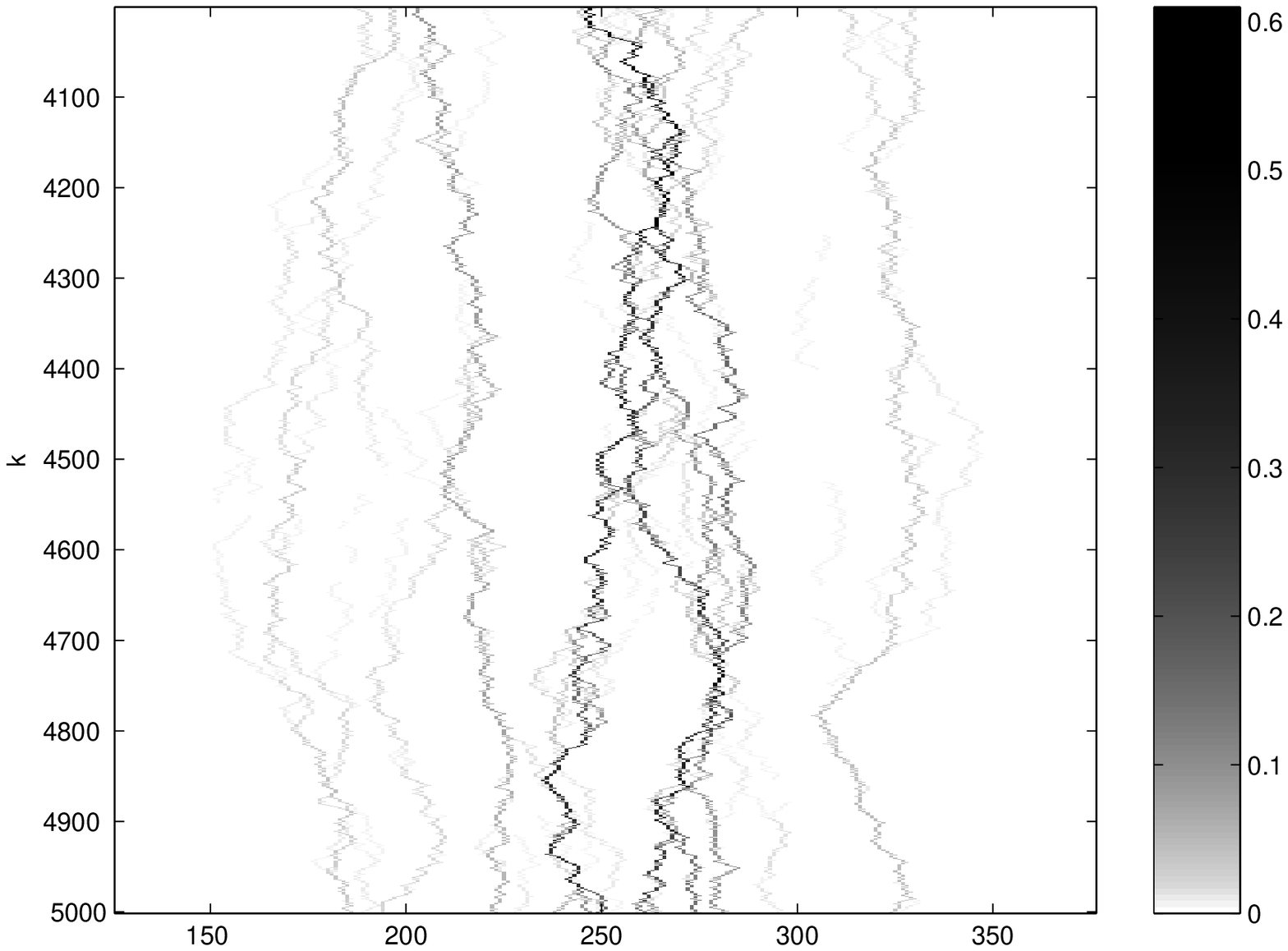}
\caption{an enlargement of the bottom  of figure \ref{fignoyau20}}}\qquad
\begin{minipage}{2.9in}
\includegraphics*[scale=0.4]{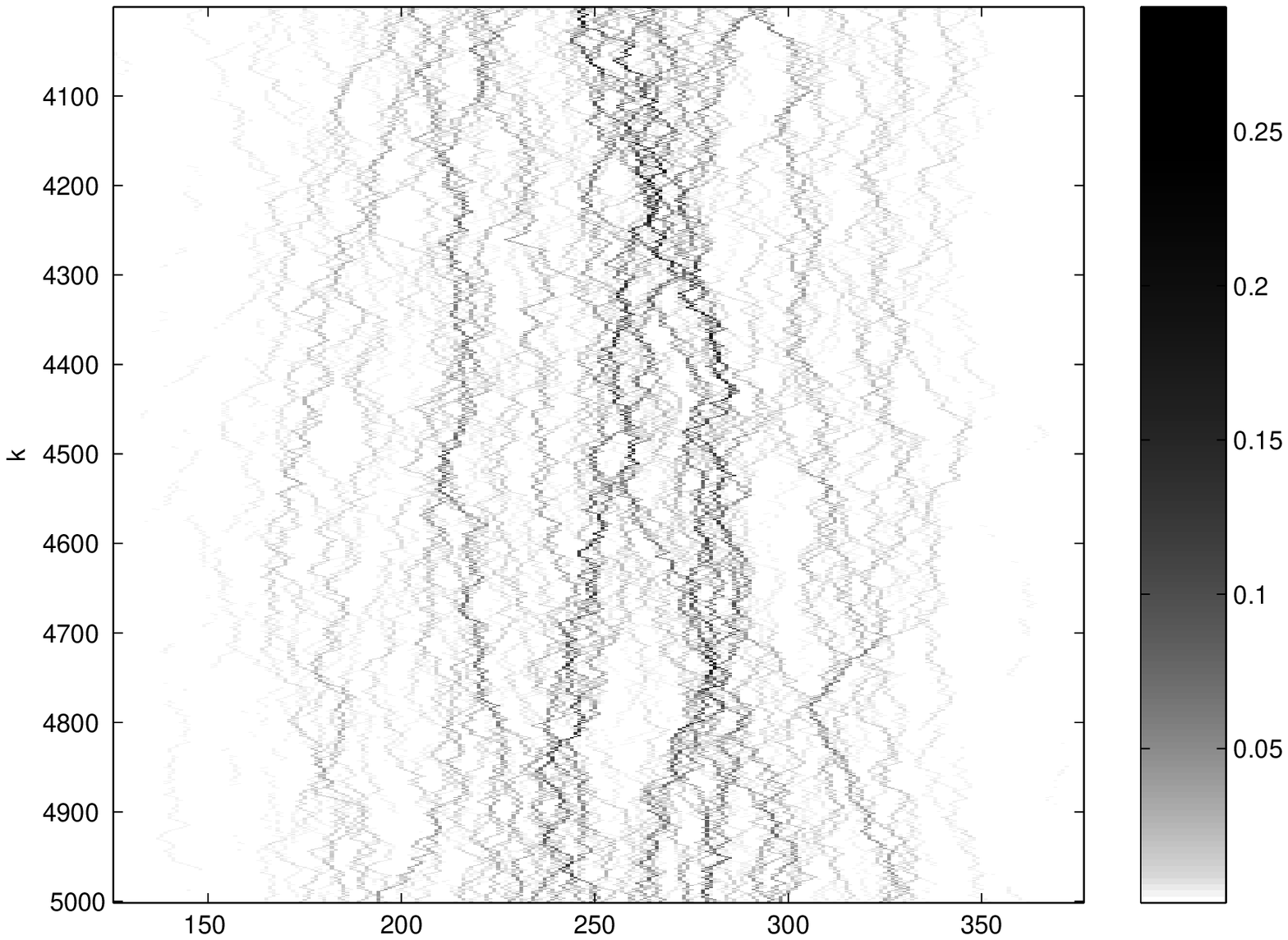}
\caption{an enlargement of the bottom of figure \ref{fignoyau100}}
\end{minipage}
\end{figure}
%
%
\subsection{Consistent system of Brownian sticky  kernels on $S^1$}
In \cite{LeJanRaimond2,LeJanRaimond3}, a general class of sticky kernels associated
with systems of coalescing particles that move as independent Levy
processes is constructed. Sticky kernels are characterized by  a parameter of
stickiness $\tau$ and the exponent $\psi$ of the associated Levy
process, that is $\psi(u)=\frac{u^2}{2}$ in the case of  Brownian
sticky kernels. Let us describe the properties of Brownian sticky kernels
on $S^1$. 

For $n\in\N^*$, let us consider the measure $m^{(n)}$ defined on
$(S^{1})^n$ by $m^{(n)}=\sum_{\pi\in\cP_n} p^{(a)}_{\pi}\lambda_\pi$ where
$\lambda_{\pi}$ is the image of the Lebesgue measure $\lambda^{\otimes
|\pi|}$ on $(S^{1})^{|\pi|}$ by the one-to-one map $\phi_{\pi}$ (that sends
$(S^{1})^{|\pi|}$ onto $(S^{1})^{n}$) and $\{p^{(a)}_{\pi},\
\pi\in\mathcal{P}_n\}$ is the partition function  defined in proposition \ref{decommesurinv}.
As this defines an exchangeable and consistent system of measures, it
follows from Kingman's representation theorem that there exists a
random measure on $S^{1}$ such that 
$m^{(n)}=E(\mu^{\otimes n})$ for every $n\in\N^*$.  In our
case, $\mu$ is the Dirichlet process\footnote{Let $\delta$ be a non null finite measure on a 
measurable 
space $(Y,\mathcal{Y})$. A Dirichlet
process  on $(Y,\mathcal{Y})$ of parameter $\delta$ is a random measure $\mu$ on $Y$ such that for 
any finite measurable partition $(B_{1},\ldots,B_{r})$ of $Y$, 
$(\mu(B_{1}),\ldots,\mu(B_{r}))$ has the Dirichlet law 
$\mathcal{D}(\delta(B_{1}),\ldots, \delta(B_{r}))$.} of parameter $a\lambda$ on $S^{1}$. 

Using this family of measures, a consistent system of Feller
semigroups denoted by $(P^{(n)}_t)$ can be defined via their Dirichlet
forms: 
\begin{prop}[Y. Le Jan and O. Raimond, \cite{LeJanRaimond2}]
For $k\in\N^{*}$, let $\mathcal{E}^{\odot k}$ be the Dirichlet form
defined on $L^{2}((S^{1})^k,\lambda^{\otimes k})$,
associated with $k$ independent  Brownian motions on
$S^{1}$. 
For every $n\in\N^*$, let  $\cE^{(n)}$ denote the Dirichlet form on
$C^{1}((S^{1})^n)$ defined as follows:

\begin{equation}
\label{formdiric}\forall f,g\in C^{1}((S^{1})^n),\ 
 \cE^{(n)}(f,g)=\sum_{\pi\in\mathcal{P}_n}p^{(a)}_\pi\mathcal{E}^{\odot
|\pi|}(f\circ \phi_\pi,\ g\circ \phi_\pi).
\end{equation}
\begin{itemize}
\item The semigroups $(P_{t}^{(n)})$, $n\in\N^*$ associated with the
Dirichlet forms $\cE^{(n)}$, $n\in \N^*$ define a consistent system of strong
Feller semigroups.
\item For every  $n\in\N^*$, the generator of $\cE^{(n)}$ denoted by
$A^{(n)}$ has the following expression on $C^{2}$ functions:
\begin{equation}
\label{generator}
A^{(n)}(f)=\frac{1}{2}\sum_{\pi\in\mathcal{P}_n}(\Delta^{(|\pi|)}(f\circ
\phi_{\pi}))\circ\phi_{\pi}^{-1}\un_{C_{\pi}}\quad 
\forall f\in C^{2}((S^{1})^n)
\end{equation}
\end{itemize}
\end{prop}
Since $\{(P_{t}^{(n)})_t,\ n\in\N^*\}$ is a compatible family of
Feller semigroups on $S^{1}$, it follows from theorem  1.1.4 in 
\cite{LeJanRaimond1}  that  it is possible to construct a stochastic flow of
kernels $(K_{s,t})$ such that $E(K_{0,t}^{\otimes n})=P_{t}^{(n)}$ for
every $t\in\R_+$ and $n\in\N^*$.  This stochastic flow  was named  sticky
flow of parameter $\tau=\frac{1}{a+1}$ and exponent $\psi(u)=
\frac{u^2}{2}$ (or Brownian sticky flow of parameter
$\tau=\frac{1}{a+1}$). \\

Given a probability measure $\nu_0$ on $S^{1}$, the sticky flow
induces a  stochastic process $(\nu_0K_{0,t})_t$ on the set $\mathcal{M}_1(S^{1})$ of probability measures
on $S^{1}$. 
\begin{prop}[Y. Le Jan and O. Raimond, \cite{LeJanRaimond2,LeJanRaimond3}]
Let $(K_{s,t})_{s\leq t}$ be a sticky Brownian flow of parameter
$\tau=\frac{1}{a+1}$. 
\begin{enumerate}
\item Let $\nu_{0}$ be a probability measure on
$S^{1}$. Then $(\nu_{0}K_{0,t})_t$ is a Feller process on
$\mathcal{M}_1(S^{1})$. It converges in law to the Dirichlet process
of parameter $a\lambda$ on $S^{1}$.
\item Let $\mu$ be a Dirichlet process
of parameter $a\lambda$ on $S^{1}$, independent of $(K_{0,t})_{t}$. Then
 $(\mu K_{0,t})_t$ is a reversible process.  
\end{enumerate} 
\end{prop}
\section{Convergence theorem}
We shall now establish the ``weak'' convergence   of the Beta flow
of parameter $a$ on $T_N$ to  the Brownian 
sticky flow  of parameter $\tau=\frac{1}{1+a}$. We keep the same
notations as those introduced  in section \ref{sectbetaflow}.
\begin{theo}
\label{cvnoyaux} 
Let $\l$ be a positive integer. For $N\in \N^*$, let $\nu_N$ be a probability
measure on $(\frac{1}{N}(\Z/N\Z))^{\l}$. Assume that $\nu_N$ converges weakly
to a probability measure $\nu$ on $(S^{1})^{\l}$ as $N$ tends to $+\infty$. 
\begin{itemize}
\item[(i)]
For every $t\in
\R_+$, the sequence of measure-valued random variables $\nu_NK_{N,0,t}^{\otimes \l}$
converges in law to $\nu K_{0,t}^{\otimes \l}$ as $N$ tends to
$+\infty$. 
\item[(ii)]For every continuous
  functions $f$ and $g$ on $(S^{1})^{\l}$ and for every ${t\in\R_{+}}$,  
$$\int P_{N,0,t}^{(\l)}(g)(x)f(x) d\nu_N(x) \mbox{ 
  converges to }\int P_{0,t}^{(\l)}(g)(x)f(x) d\nu(x)$$ as $N$ tends to 
$+\infty$.
\end{itemize}For $N\in \N^*$, let $\eta_N$ be a probability
measure on $\frac{1}{N}(\Z/N\Z)$. Assume that $(\eta_N)_N$ converges weakly
to a probability measure $\eta$ on $S^{1}$. 
\begin{itemize}
\item[(iii)] The finite dimensional distributions of the
  measure-valued Markov process \linebreak[4] ${\{\eta_NK_{N,0,t},\ t\in\R_{+}\}}$ weakly converge to the
finite  dimensional distributions of \linebreak[4]${\{\eta K_{0,t},\ t\in\R_{+}\}}$ as $N$ tends to $+\infty$.
\end{itemize}
\end{theo}   
\rems 
\begin{enumerate}
\item Assertions of the theorem also hold  with  sequences of
probability measures $\nu_N$ and $\eta_N$ defined on
$(\frac{1}{2N}((2\Z+1)/2N\Z))^{\l}$ and  $\frac{1}{2N}((2\Z+1)/2N\Z)$
respectively that converge weakly as $N$ tends to $+\infty$. 
\item  
 For $N\in \N^*$, let $\nu^{(1)}_N,\ldots,\nu^{(r)}_N$  be  probability
measures on $\frac{1}{N}(\Z/N\Z)$ such that $(\nu^{(i)}_N)_N$ converges weakly
to a probability measure $\nu^{(i)}$ of $S^{1}$ for every ${i\in\{1,\ldots,r\}}$. It follows from assertion $(i)$ of
theorem \ref{cvnoyaux} that 
$(\nu^{(1)}_NK_{N,0,t},\ldots,\nu^{(r)}_NK_{N,0,t})$ converges in law to
$(\nu^{(1)}K_{0,t},\ldots,\nu^{(r)}K_{0,t})$ for every
$t\in\R_{+}$. This property could be used as a definition of the
convergence in law of stochastic flows of kernels with independent
increments.   
\end{enumerate}

Before going into details, let us explain the scheme of the 
proof of  theorem \ref{cvnoyaux}.   First,  $(ii)$ and $(iii)$ will be deduced from assertion
$(i)$ of the theorem. 
The main step of the proof of $(i)$ is to show the following
convergence 
of the resolvents:
\begin{prop}
\label{cvresw}Let $n\in\N^*$.  
Let $(V_{N,\alpha}^{(n)})_{\alpha>0}$ denote the resolvent
associated with  the
$n$-point motion of the Beta flow of parameter $a$ on $T_{N}^{(n)}$ and let $(V_{\alpha}^{(n)})_{\alpha>0}$
denote the resolvent of the $n$-point motion of the Brownian sticky flow of
parameter $\tau=\frac{1}{a+1}$.   
 For every $\alpha\in \R_{+}^*$ and  continuous functions $f$ and $g$ on $(S^{1})^{n}$, 
$\displaystyle{\int V_{N,\alpha}^{(n)}(f)gdm_{N}^{(n)}}$ converges to $\displaystyle{\int
V_{\alpha}^{(n)}(f)gdm^{(n)}}$ as $N$ tends to $+\infty$.
\end{prop} 
As the discrete and the continuous $n$-point processes are both
reversible, an argument using spectral measures allows to deduce the
following weak convergence of the semigroups from the weak convergence of the
resolvents:
\begin{prop}
\label{cvmoment}    
Let $n\in\N^*$. If  $f$ and $g$ are  continuous functions on
$(S^{1})^n$, then $\int gP^{(n)}_{N,t}(f)dm_{N}^{(n)}$ converges to 
    $\int gP^{(n)}_{t}(f)dm^{(n)}$ as $N$ tends to $+\infty$. 
\end{prop}
A last step consists in making the previous convergence result of the
semigroups also  true if $m_{N}^{(n)}$ is replaced by
 any sequence of probability measures $\nu_N$ on
 $(\frac{1}{N}(\Z/N\Z))^{n}$ that weakly converges.
\subsection{Convergence of the resolvents of the $n$-point motions}
 The convergence of the resolvents  is based on the convergence of the invariant measures
$m_{N}^{(n)}$ and of the generator of the $n$-point motion of the 
flow of Beta matrices together with  a Lipschitz property of  the discrete resolvent 
$V_{N,\alpha}^{(n)}$.
 Before proving proposition \ref{cvresw}, let us give precise
statements of these three points.
\subsubsection{Convergence of the invariant measures}
\begin{lem}
    \label{convmes}
For every partition $\pi$ of $[n]$, let $f_{\pi}$ be a function
defined on $(S^{1})^n$,  Lipschitz on $E_\pi$ and  vanishing outside 
$E_\pi$. Let $f=\sum_{\pi\in\mathcal{P}_n}f_{\pi}$.  Then $$\Big|\int f dm^{(n)}_N-\int
fdm^{(n)}\Big|\leq \frac{C_n}{N}\sum_{\pi\in\mathcal{P}_n}(||f_{\pi}||_{Lip}+||f_{\pi}||_{\infty}).$$ 
\end{lem} 
\begin{proof} Let $E_{N,\pi}$ denote the intersection of $E_{\pi}$ 
with $T_{N}^{(n)}$ and let $\lambda_{N,\pi}$ be the uniform 
distribution on $E_{N,\pi}$. By formula (\ref{mesinvbeta}) and
proposition \ref{decommesurinv}, ${m^{(n)}_N=\sum_{\pi\in
  \mathcal{P}_n}p^{(a)}_{\pi}\lambda_{N,\pi}}$. As the measure $m^{(n)}$ has a
similar decomposition:  ${m^{(n)}=\sum_{\pi\in
  \mathcal{P}_n}p^{(a)}_{\pi}\lambda_{\pi}}$, it  suffices to prove that for all partitions $\pi,\pi'$ of $[n]$, 
$$\Big|\int_{E_{N,\pi'}} f_{\pi} d\lambda_{N,\pi'}-\int_{E_{\pi'}}
f_{\pi}d\lambda_{\pi'}\Big|\leq
\frac{C_n}{N}(||f_{\pi}||_{Lip}+||f_{\pi}||_{\infty}).$$
Let us first consider the case  $\pi$ is not finer than $\pi'$,
that is there is a nonempty block $B$ of $\pi$ intersecting at least
two blocks of $\pi'$. Let $\tilde{\pi}$ be the partition of $[n]$
obtained by merging   the blocks of $\pi$ that intersect the same block of
$\pi'$. Then  $E_{\pi'}\cap E_{\pi}$ is a subset
of $E_{\tilde{\pi}}$. As $\tilde{\pi}$ is a coarser partition than
$\pi'$, $\lambda_{\pi'}(E_{\tilde{\pi}})=0$ and $|E_{N,\tilde{\pi}}|\leq
\frac{1}{N}|E_{N,\pi'}|$. Thus 
$$\Big|\int_{E_{N,\pi'}}f_{\pi} d\lambda_{N,\pi'}-\int_{E_{\pi'}}
f_{\pi}d\lambda_{\pi'}\Big|=\frac{1}{|E_{N,\pi'}|}\Big|\sum_{x\in 
E_{N,\pi'}\cap E_{N,\pi}}f_{\pi}(x)\Big|\leq \frac{1}{N}||f_{\pi}||_{\infty}.$$
Let us now consider the case $\pi$ is equal to $\pi'$ or finer than $\pi'$. Then
$E_{\pi'}\subset E_{\pi}$. If $\pi'$ has $k$ nonempty
blocks then  
\begin{multline*}\Big|\int_{E_{N,\pi'}} f_{\pi} d\lambda_{N,\pi'}-\int_{E_{\pi'}}
f_{\pi}d\lambda_{\pi'}\Big|=\Big|\frac{1}{2N^k}\sum_{x\in
T^{(k)}_{N}}f_{\pi}(\phi_{\pi'}(x))-\int_{(S^{1})^k}f_{\pi'}(\phi_{\pi'}(x))dx\Big|.
\end{multline*}
The function $f_{\pi'}\circ \phi_{\pi'}$ is a Lipschitz function with Lipschitz
coefficient smaller than $n||f_{\pi'}||_{Lip}$. Thus it remains to
establish the following result: for every  $k\in \N^*$, there exists a constant $C_k$ such that for every
Lipschitz function $g$ on $(S^{1})^k$, 
$$\Big|\frac{1}{2N^k}\sum_{x\in T^{(k)}_N}g(x)-\int_{(S^{1})^k}g(x)dx\Big|\leq
\frac{C_k}{N}||g||_{Lip}.$$
The proof can be made by induction on $k$.
 \end{proof}
\subsubsection{Convergence of the generators}
\begin{lem}
    \label{convgen}
For every $n\in\N^*$, let  $A^{(n)}_{N}$  denote  the generator of the $n$-point motion of 
the Beta flow on $T^{(n)}_N$ of parameter $a$. For every $C^2$ function $f$ on
$(S^1)^n$,
$$\sup_{x\in
T^{(n)}_N}|A^{(n)}_{N}(f)(x)-A^{(n)}(f)(x)|$$
converges to $0$ as $N$ tends to $+\infty$. 
\end{lem}
\begin{proof}  
Let $n\in\N^*$ and $f$ be a $C^{2}$ function defined on $(S^1)^n$. Let
us recall the expression of $A^{(n)}(f)$ :
$$A^{(n)}(f)=\frac{1}{2}\sum_{\pi\in\mathcal{P}_n}
\Delta_{\pi}(f)\un_{C_\pi}\mbox{  where }
\Delta_{\pi}(f)(x)=\Delta^{(|\pi|)}(f\circ
\phi_{\pi})(\phi_{\pi}^{-1}(x)).$$
Thus it suffices to prove that for every partition  $\pi$ of  $[n]$,
    $$\displaystyle{\sup_{x\in C_{\pi}\cap T_{N}^{(n)}}\mid
2A_{N}^{(n)}(f)(x) -\Delta_{\pi}f(x)\mid}$$ converges
to $0$  as $N$ tends to $+\infty$.\\
We shall introduce the discrete version of 
the Laplacian $\Delta_{\pi}$, denoted by $\Delta_{N,\pi}$ and we shall prove that the difference 
between $2A_{N}^{(n)}(f)$ and  $\Delta_{N,\pi}f$ goes to zero as $N$ 
tends to $+\infty$. \\
The expression linking  the generator $A^{(n)}_N$ and  the transition
matrix $P_{N}^{(n)}$ is the following:
 for a function  $g$ defined on  $(S^{1})^n$ and $x\in T^{(n)}_N$,  
\[ A^{(n)}_{N}(g)(x)=4N^{2}\left(\sum_{\varepsilon\in\{\pm 
\frac{1}{2N}\}^n}P^{(n)}_{N}(x,x+\varepsilon)g(x+\varepsilon)-
g(x)\right).\]
Let us recall the expression of $P^{(n)}_{N}$: for every $x\in
T_{N}^{(n)}$  and $\eps\in\{\pm \frac{1}{2N}\}^{n}$, 
\[ P^{( n )}_{N} ( x, x+\varepsilon )  = \prod_{l \in \Z/2N\Z}
   \frac{\prod_{i = 0}^{s_l^+ ( x, \varepsilon ) - 1} (
   \frac{a}{2N} + i ) \prod_{i = 0}^{s_l^- (x, \varepsilon
   ) - 1} ( \frac{a}{2N} + i )}{\prod_{i = 0}^{s_l ( x ) -
   1} ( \frac{a}{N} + i )} . \]
where $s_l ( x )$ denotes the number of coordinates of $2Nx$ equal to $l$
   and $$s_l^{\pm} ( x, \varepsilon ) = \Card ( \{ i \in \{ 1, \ldots, n \}, 2Nx_i
= l \text{ and } 2N\varepsilon_i = \pm 1 \} ).$$
As $s^{+}_{l}(x,\varepsilon)=s^{-}_{l}(x,-\varepsilon)$, 
$P^{(n)}_{N}(x,x+\varepsilon)=P^{(n)}_{N}(x,x-\varepsilon)$. Thus, 
\[ 2A^{(n)}_{N}(g)(x)=\sum_{\varepsilon\in\{\pm 
\frac{1}{2N}\}^n}P^{(n)}_{N}(x,x+\eps)L_{n}(g)(x,\eps)\ 
\forall x\in T^{(n)}_N,\]
  where $L_{n}(g)(x,h)=4N^{2}(g(x+h)+g(x-h)-2g(x))$ for $h\in \R$.\\
 For $d,N\in \N^{*}$, let $\Delta_{N}^{(d)}$ denote the discrete Laplacian on 
 the set $T^{(d)}_{N}$. Its expression in terms of $L_{d}$ is:  
 $$\Delta_{N}^{(d)}g=\frac{1}{2^d}\sum_{\epsilon\in\{\pm 
 1\}^d}L_{d}(g)(\ \cdot\ ,\frac{\epsilon}{2N})$$ for a 
    function $g$ defined on $T^{(d)}_{N}$. 
 Let $\pi$ be a partition of $[n]$ having $d$ nonempty blocks and set 
  $\Delta_{N,\pi}g(\cdot)=\Delta_{N}^{( d)}(g\circ \phi_{\pi})
  (\phi^{-1}_{\pi}(\cdot))$ for a function $g$ defined on $T^{(n)}_N$.
As the restriction of $\phi_{\pi}$ to $\{\pm 1\}^{d}$ is a one-to-one
map onto ${E_{\pi}\cap\{\pm 1\}^{n}}$, 
$$\Delta_{N,\pi}f(x)=\frac{1}{2^{d}}\sum_{\varepsilon\in\{\pm
\frac{1}{2N}\}^n\cap E_{\pi}}L_{n}(f)(x,\varepsilon).$$ 
Consequently, the expression of 
$2A_{N}^{(n)}f(x) -\Delta_{\pi}f(x)$ can be split
into the three following terms:
\begin{eqnarray*}
I^{1}_N(x)&=& \sum_{\varepsilon \in E_\pi\cap \{\pm
\frac{1}{2N}\}^n}(P^{(n)}_{N}(x,x+\varepsilon)-\frac{1}{2^{d}})L_{n}(f)(x,\varepsilon)\\
I^{2}_N(x)&=& \Delta_{N,\pi}f(x)-\Delta_{\pi}f(x)\\
I^{3}_N(x)&=& \sum_{\varepsilon \in E^{c}_\pi\cap \{\pm
\frac{1}{2N}\}^n}P^{(n)}_{N}(x,x+\eps)L_{n}(f)(x,\varepsilon)
\end{eqnarray*}
\paragraph{Asymptotic behaviour of $I^{2}_{N}(x)$.}
For a $C^{2}$ function $g$ on
$(S^1)^r$,  $x\in(S^{1})^{r}$ and $\epsilon\in\{\pm 1\}^{r}$,  set
$\displaystyle{\mathcal{L}_r(g)(x,\epsilon)=\sum_{i=1}^{r}
\sum_{j=1}^{r}\epsilon_{i}\epsilon_{j}\partial^{2}_{i,j}g(x)}.$
 It follows from a Taylor expansion with integral remainder that   
 ${|L_{r}(g)(x,\frac{\epsilon}{2N})-\mathcal{L}_r(g)(x,\epsilon)|}$ converges to
zero uniformly on $x\in T^{(r)}_N$ and $\epsilon \in\{\pm 1\}^r$ as
 $N$ tends to $+\infty$. Thus 
$I^{2}_{N}(x)=\Delta_{N,\pi}f(x)-\Delta_{\pi}f(x)$
 converges to zero uniformly on $x\in T^{(n)}_{N}$.
\paragraph{Asymptotic behaviour of $P^{(n)}_{N}(x,x+\varepsilon)$.}
 Let $(B_1,\ldots,B_d)$ denote the
non\-empty blocks of $\pi$ and let  $x\in C_{\pi}\cap T^{(n)}_N$. Then 
$$P^{(n)}_{N}(x,x+\varepsilon)=\prod_{l=1}^{d}\frac{\gamma_{\frac{a}{2N}}(
\underset{i\in
B_l}{\sum}\tmmathbf{1}_{\{\eps_i=\frac{1}{2N}\}})\gamma_{\frac{a}{2N}}
(\underset{i\in
B_l}{\sum}\tmmathbf{1}_{\{\eps_i=-\frac{1}{2N}\}})}{\gamma_{\frac{a}{N}}(|B_l|)}\ \forall \eps\in\{\pm \frac{1}{2N}\}^n$$ where $\gamma_b(u)$
denotes the product 
$\prod_{i=0}^{u-1}(b+i)$ for $b>0$ and $u\in \N$ with the convention $\prod_{i=0}^{-1}=1$.
A computation shows that 
\begin{itemize}
\item if $u,v\in \N^*$ then 
$$\frac{\gamma_{\frac{a}{2N}}(u)\gamma_{\frac{a}{2N}}(v)}{\gamma_{\frac{a}{N}}(u+v)}\leq
\frac{a}{4N}\frac{u!v!}{(u+v-1)!} \mbox{ for }2N\geq a.$$ 
\item  if $u\in\N^*$ and $v=0$ then 
\begin{eqnarray*}|\frac{\gamma_{\frac{a}{2N}}(u)\gamma_{\frac{a}{2N}}(0)}{\gamma_{\frac{a}{N}}(u)}-\frac{1}{2}|
&=&\frac{1}{2}(1-\prod_{i=1}^{u-1}(1-\frac{a}{2a+2Ni}))\\
&\leq&
\frac{1}{2}(1-(1-\frac{a}{2a+2N})^{u-1}).
\end{eqnarray*}
\end{itemize}
Thus
$$\sup_{x\in C_{\pi},\ \eps\in\{\pm \frac{1}{2N}\}^{n}\cap
E^{c}_{\pi}}P^{(n)}_{N}(x,x+\varepsilon) \mbox{ and} \sup_{x\in C_{\pi},\
\eps\in\{\pm \frac{1}{2N}\}^{n}\cap
E_{\pi}}|P^{(n)}_{N}(x,x+\varepsilon)-\frac{1}{2^d}|$$ converge to $0$ as $N$ tends to
$+\infty$.  
\end{proof}    
\subsubsection{Lipschitz property of the resolvents and the 
semigroups associated with the
Beta flow}
\begin{lem}\label{propresolv}Let $(P^{(n)}_{N,t})_t$ denote the $n$-point Markovian semigroup of 
the Beta flow on $T^{(n)}_N$. \\
If $f: (S^{1})^n\rightarrow \mathbb{R}$ is  a Lipschitz function
then
\begin{itemize}
\item the Lipschitz 
coefficient of $P_{N,t}^{(n)}(f)$ is  bounded by $|| f||_{Lip}$.
\item the Lipschitz coefficient of 
$V^{(n)}_{N,\alpha}(f)$ is 
 bounded by $\frac{1}{\alpha}|| f||_{Lip}$.
\end{itemize}
\end{lem}
\begin{proof}
We use a coupling argument borrowed from \cite{LeJanRaimond2}. Let
$\underline{x}\!=\!(x_1,\ldots,x_{n+1})$ be  a point of $T^{(n+1)}_{N}$ such that
$x_1\neq x_2$ and   let
$X_{t}=(X_{t}^{(1)},\ldots,X_{t}^{(n+1)})$ be a Markov chain on
$T^{(n+1)}_{N}$ with transition matrix $P_{N}^{(n+1)}$ and with
initial 
point $\underline{x}$. Set $\tau=\inf\{s>0,\
X^{(1)}_{s}=X^{(2)}_{s}\}$. Since $(X^{(1)}_{t},X^{(2)}_{t})$ is a positive 
recurrent 
Markov chain on $T^{(2)}_{N}$, $\tau$ is almost surely finite. Let us define
two processes $(Y_t)_t$ and $(Z_{t})_t$  on $T^{(n)}_N$:
\begin{itemize}
\item $Y_{t}^{(1)}=X^{(1)}_t$ and  $Y_{t}^{(i)}=X_{t}^{(i+1)}$ for
$i\in\{2,\ldots,n\}$,
\item $Z_{t}^{(1)}=X^{(2)}_{t}\un_{t\leq \tau}+X^{(1)}_{t}\un_{t> \tau}$ and $Z_{t}^{(i)}=X_{t}^{(i+1)}$ for
$i\in\{2,\ldots,n\}$.
\end{itemize}
As  $\{(P^{(n)}_{N,t})_t,\ n\in\N^*\}$ defines a consistent family of
Markovian semigroups, $(Y_t)_t$ and 
${(X^{(2)}_t,X^{(3)}_t,\ldots,X^{(n+1)}_t)_t}$ are both Markov  processes
with semigroup $(P^{(n)}_{N,t})_t$. The strong Markov property implies
that $(Z_t)_t$ is also a Markov process with semigroup
$(P^{(n)}_{N,t})_t$. As $Y_t=Z_t$ if $t\geq \tau$,
 \begin{multline*}
|P^{(n)}_{N,t}(f)(x_1,x_3,\ldots,x_{n+1})-
P^{(n)}_{N,t}(f)(x_2,x_3,\ldots,x_{n+1})|\\=|E(f(Y_{t\wedge
\tau})-f(Z_{t\wedge \tau}))|
\leq
||f||_{Lip}E(d(X^{(1)}_{t\wedge \tau}, X^{(2)}_{t\wedge \tau})).
\end{multline*} and 
\begin{multline*}
|V^{(n)}_{N,\alpha}(f)(x_1,x_3,\ldots,x_{n+1})-
V^{(n)}_{N,\alpha}(f)(x_2,x_3,\ldots,x_{n+1})|\\=|\int_{0}^{+\infty}E(f(Y_{t\wedge
\tau})-f(Z_{t\wedge \tau}))e^{-\alpha t}dt|\\
\leq
||f||_{Lip}\int_{0}^{+\infty}e^{-\alpha t}E(d(X^{(1)}_{t\wedge \tau}, X^{(2)}_{t\wedge \tau}))dt.
\end{multline*}
Let us show that for every $t\geq 0$, $E(d(X^{(1)}_{t\wedge \tau}, X^{(2)}_{t\wedge \tau}))
\leq d(x_1,x_2)$.
Without loss of generality, one may assume that
$w=x_1-x_2\in\{0,\frac{1}{N},\ldots,\frac{N-1}{N}\}$. Let $(\hat{X}^{(1)}_{t},\hat{X}^{(2)}_{t})$
be a Markov chain on $(\frac{1}{2N}\Z)^2$ starting from $(x_1,x_2)$ whose transition
matrix $\hat{P}$ is defined by:\\
\centerline{$\hat{P}(x+k,y+l)=P^{(2)}_N(x,y)$ for every 
$(x,y)\in T^{(2)}_N$ and $k,l\in\Z^2$.} Set $W_t=\hat{X}^{(1)}_t-\hat{X}^{(2)}_t$.
 Since $2Nw$ is
even,  for
every $t\geq 0$, $W_{t\wedge \tau}$ remains nonnegative, whence  ${E(d(\hat{X}^{(1)}_{t\wedge\tau},\hat{X}^{(2)}_{t\wedge \tau}))\leq E(W_{t\wedge \tau})}$. As
$(W_t)$ is a martingale, for every $t\geq 0$, ${E(W_{t\wedge
  \tau})=w}$. 
Therefore, for every $t\geq 0$, $$E(d(X^{(1)}_{t\wedge \tau}, X^{(2)}_{t\wedge \tau}))=E(d(\hat{X}^{(1)}_{t\wedge
\tau},\hat{X}^{(2)}_{t\wedge \tau}))\leq w=d(x_1,x_2).$$
 We have obtained the following inequalities: for every $t\geq 0$,
\begin{alignat*}{2}|P^{(n)}_{N,t}(f)(x_1,x_3,\ldots,x_{n+1})-
	P^{(n)}_{N,t}(f)(x_2,x_3,\ldots,x_{n+1})|\!&\leq\!
||f||_{Lip}d(x_1,x_2),\\ 
|V^{(n)}_{N,\alpha}(f)(x_1,x_3,\ldots,x_{n+1})-
V^{(n)}_{N,\alpha}(f)(x_2,x_3,\ldots,x_{n+1})|\!&\leq\!
\frac{1}{\alpha}||f||_{Lip}d(x_1,x_2).
\end{alignat*}As the semigroup
$(P^{(n)}_{N,t})_t$ is invariant
by the action of a permutation, for every 
$x,y\in T_{N}^{(n)}$, 
\begin{eqnarray*}&& |P^{(n)}_{N,t}(f)(x)-P^{(n)}_{N,t}(f)(y)|\leq
||f||_{Lip}\sum_{i=1}^{n}d(x_i,y_i),\\
&&|V^{(n)}_{N,\alpha}(f)(x)-V^{(n)}_{N,\alpha}(f)(y)|\leq
\frac{1}{\alpha}||f||_{Lip}\sum_{i=1}^{n}d(x_i,y_i).\quad  
\end{eqnarray*} 
 \end{proof}
\subsubsection{Proof of  proposition \ref{cvresw}. }
A density argument reduces the problem to showing that  for all 
$C^{1}$ functions $f$ and $g$ on $(S^{1})^{n}$, 
$\int V_{N,\alpha}^{(n)}(f)gdm_{N}^{(n)}$ converges to $\int
V_{\alpha}^{(n)}(f)gdm^{(n)}$ as $N$ tends to $+\infty$. \\
Let us introduce an extension of $V_{N,\alpha}^{(n)}(f)$ to 
$(S^1)^n$:  
\begin{lem}
\label{extfunct}
A  function $g$ on $T_{N}^{(n)}$ can be extended to a Lipschitz function 
$\tilde{g}$ on $(S^{1})^{n}$ such that: \begin{itemize}
\item $\mid \mid \tilde{g} \mid \mid_{\infty}=\mid \mid g \mid \mid_{\infty}$ 
\item $\mid \mid \tilde{g} \mid \mid_{Lip}\leq C_{n}\mid \mid g \mid 
\mid_{ Lip}$ where $C_n$ is a constant only depending on $n$. 
\item $\tilde{g}$ is differentiable on $(S^1)^{n}-
\mathcal{R}$ where $\mathcal{R}$ is the subset  of points having
at least one coordinate in  $\frac{1}{2N}(\mathbb{Z}/2N\mathbb{Z})$ and $||\partial_i
\tilde{g}(x)||\leq C_n ||g||_{Lip}$ for all $i\in\{1,\ldots,n\}$ and
$x\in (S^1)^{n}-\mathcal{R}$. 
\end{itemize}
\end{lem}
\begin{proof}
Firstly, a function $g$ on  $T_{N}^{(n)}$ is extended to a
function $\bar{g}$ on the 
lattice $(\frac{1}{2N}(\Z/2N\Z))^n$ as follows: for $x\in
(\frac{1}{2N}(\Z/2N\Z))^n-T_{N}^{(n)}$, set
$$\bar{g}(x)=\frac{1}{|V_{x}|}\sum_{y\in V_x}g(y)$$ where $V_{x}$ is the
set of the nearest points of $x$ in $(\frac{1}{2N}(\Z/2N\Z))^n$
in the sense of the  distance
$d_n(x,y)=\sum_{i=1}^{n}d(x_i,y_i)$. This
extension has the following properties: 
\begin{itemize}
\item $||\bar{g}||_{\infty}=||g||_{\infty}$,
\item there is a constant $C_n>0$ such that for every  function
$g:T^{(n)}_N\rightarrow \R$, ${||\bar{g}||_{Lip}\leq
C_n||g||_{Lip}}$. 
\end{itemize}
Lastly, a function $f$ defined on $(\frac{1}{2N}(\Z/2N\Z))^n$ is
extended to a function $\hat{f}$ on  $(S^{1})^n$ as follows. A point 
$x=(x_1,\ldots,x_n)$ in an elementary cube
$\prod_{i=1}^{n}]\frac{k_i}{2N},\frac{k_i+1}{2N}[$ is the barycentre of
the vertices of this cube
$\frac{k+\eta}{2N}$, $\eta\in\{0,1\}^n$ with the weights
$$\alpha_{n}(k+\eta,x)=(2N)^{n}\prod_{i=1}^{n}(x_i-\frac{k_i}{2N})^{\eta_i}
(\frac{k_i+1}{2N}-x_i)^{1-\eta_i}$$
respectively. Then we set $\hat{f}(x)$ as the convex combination of
the points $f(\frac{k+\eta}{2N})$ with the
weights $\alpha_{n}(k+\eta,x)$ for every
$\eta\in\{0,1\}^n$ :
$$\hat{f}(x)=\sum_{\eta\in\{0,1\}^n}\alpha_{n}(k+\eta,x)f(\frac{k+\eta}{2N}).$$
Let us list some properties of this extension: 
 \begin{itemize}
\item $||\hat{f}||_{\infty}=||f||_{\infty}$ and $f$ is differentiable in 
$(S^{1})^{n}- \mathcal{R}$.
\item   $||\hat{f}||_{Lip}\leq
||f||_{Lip}$ and  for every $i\in\{1,\ldots,n\}$,
$|\partial_i\hat{f}(x)|\leq ||f||_{Lip}$ if $x\in (S^{1})^{n}- \mathcal{R}$. 
\end{itemize}
 \end{proof}
By lemma \ref{propresolv}, the Lipschitz coefficient of $V^{(n)}_{N,\alpha}(f)$ 
is    bounded by a constant irrespective of $N$, hence it is also the case for 
$\widetilde{V^{(n)}_{N,\alpha}(f)}$. Thus on
 applying lemma \ref{convmes}, we obtain that the 
difference between  $\int V_{N,\alpha}^{(n)}(f)gdm_{N}^{(n)}$ 
and $\int \widetilde{V_{N,\alpha}^{(n)}(f)}gdm^{(n)}$ converges 
to zero as $N$ tends to $+\infty$. \\
The remainder of the proof is the subject of the following key lemma:
\begin{lem}
\label{propext}Let $\alpha>0$.
If $f$ is a $C^{1}$ function  on $(S^{1})^n$ then for every $g \in L^{2}(m^{(n)})$,
$\int g\widetilde{V^{(n)}_{N,\alpha}(f)}dm^{(n)}$  converges to 
$\int gV^{(n)}_{\alpha}(f)dm^{(n)}$.      
\end{lem}
\begin{proof}
For a positive real $\alpha$, set
$\mathcal{E}^{(n)}_{\alpha}(\ \cdot,\ \cdot)=\alpha\langle\ \cdot ,\cdot\
\rangle_{m^{(n)}}+\mathcal{E}^{(n)}(\ \cdot,\ \cdot)$ and $\mathcal{H}^{(n)}$
the closure of $C^{1}((S^{1})^{n})$ for the metric
$\mathcal{E}^{(n)}_1$.\\ As $\langle
u,v\rangle_{m^{(n)}}=\mathcal{E}^{(n)}_{\alpha}(V^{(n)}_{\alpha}(u),v)$
for every $u\in L^{2}(m^{(n)})$ and $v\in\mathcal{H}^{(n)}$, 
the weak convergence in the Dirichlet space 
$(\mathcal{H}^{(n)},\mathcal{E}^{(n)}_{\alpha})$ implies the weak
convergence
 in $(L^{2}(m^{(n)}),\langle\ \cdot ,\cdot\
\rangle_{m^{(n)}})$. 
Thus, it suffices to prove that for all $g\in 
\mathcal{H}^{(n)}$,  
$\mathcal{E}^{(n)}_{\alpha}(\widetilde{V^{(n)}_{\alpha,N}(f)},g)$  tends to
$\mathcal{E}^{(n)}_{\alpha}(V^{(n)}_{\alpha}(f),g)$. \\
 It
follows from 
lemmas \ref{propresolv} and \ref{extfunct} that
$\widetilde{V^{(n)}_{\alpha,N}(f)}$ is bounded by
$\frac{||f||_{\infty}}{\alpha}$  and differentiable on
$(S^{1})^{n}-\mathcal{R}$ with partial derivatives bounded by
$||f||_{Lip}$. Thus expression (\ref{formdiric}) of $\mathcal{E}^{(n)}$
enables us to establish the following inequality
$$\mathcal{E}_{\alpha}^{(n)}(\widetilde{V^{(n)}_{\alpha,N}(f)})\leq
C_{n,\alpha}(||f||^{2}_{Lip}+||f||^{2}_{\infty})$$ where
$C_{n,\alpha}$ is a constant irrespective of 
 $N$ and $f$. 
 As the set of $C^{3}$ functions is dense in 
$\mathcal{H}^{(n)}$, this 
reduces the problem to proving the convergence of $\mathcal{E}^{(n)}_{\alpha}(\widetilde{V^{(n)}_{\alpha,N}(f)},g)$ for any $C^{3}$ function 
$g$.\\
Let $g$ be  a $C^{3}$ function. The difference 
$\mathcal{E}^{(n)}_{\alpha}(\widetilde{V^{(n)}_{N,\alpha}(f)},g)-
\mathcal{E}^{(n)}_{\alpha}(V^{(n)}_{\alpha}(f),g)$ is the sum of two terms: 
\begin{eqnarray*}
    I_{N}^{(1)}&=&\mathcal{E}^{(n)}_{\alpha}(\widetilde{V^{(n)}_{N,\alpha}(f)},g)
    -\mathcal{E}^{(n)}_{N,\alpha}(V^{(n)}_{N,\alpha}(f),g)\\
    I_{N}^{(2)}&=&\int f(x)g(x)dm^{(n)}_{N}(x)-\int f(x)g(x)dm^{(n)}(x)
\end{eqnarray*}
By lemma \ref{convmes}, $I^{(2)}_{N}$ goes to zero as $N$ tends to 
$+\infty$. \\ 
Let us split up  $I_{N}^{(1)}$:
\begin{alignat*}{2}
    I^{(1)}_{N}&=&\int \widetilde{V^{(n)}_{N,\alpha}(f)}(\alpha -A^{(n)})(g) 
    dm^{(n)}-
    \int \widetilde{V^{(n)}_{N,\alpha}(f)}(\alpha -A^{(n)})(g) dm^{(n)}_{N}\\
    && +\int V^{(n)}_{N,\alpha}(f)(A^{(n)}_{N}-A^{(n)})(g) dm_{N}^{(n)}\\
\end{alignat*}   
 It follows from  expression (\ref{generator})
of the generator  $A^{(n)}$ that  $A^{(n)}(g)$ is a sum of 
Lipschitz functions on $E_{\pi}$ that vanish out of $E_{\pi}$. Thus 
 lemma \ref{convmes} can be applied to
 $\widetilde{V^{(n)}_{N,\alpha}(f)}(\alpha -A^{(n)})$. Finally, the last integral is bounded 
by $$\frac{\Arrowvert f\Arrowvert_{\infty}}{\alpha}
\sup_{x\in T_{N}^{(n)}}\mid (A^{(n)}_{N} 
-A^{(n)})(g)(x)\mid.$$ Thus it converges to $0$ as $N$ tends to
$+\infty$, by lemma \ref{convgen}.
\end{proof}
\subsection{Proof of  proposition \ref{cvmoment}.}
 Let $f$ be a continuous function on $(S^{1})^n$.
    As  the two $n$-point 
    processes are reversible Markov processes, their generators 
    $A_{N}^{(n)}$ and $A^{(n)}$  are  self-adjoint operators  on the Hilbert spaces 
    $L^{2}(m^{(n)}_{N})$ and $L^{2}(m^{(n)})$ with nonpositive spectra. Let $\nu_{N}^{f}$ and 
    $\nu^{f}$ denote the spectral measures of  $A_{N}^{(n)}$ and $A^{(n)}$ 
    respectively associated with the function $f:$ $$ \langle f,\ 
    \psi(A_{N}^{(n)})f\rangle_{m^{(n)}_N} 
    =\int_{\R_{-}}\psi d\nu^{f}_{N} \mbox{ and } \langle f,\ 
    \psi(A^{(n)})f\rangle_{m^{(n)}} 
    =\int_{\R_{-}}\psi d\nu^{f}$$ for every 
     continuous function $\psi$ on  $\R_{-}$.\\ The 
     relation between the resolvent and the generator, given in the
    discrete case by   
     ${V^{(n)}_{N,\alpha}=(\alpha -A_{N}^{(n)})^{-1}}$, and the 
     convergence of the resolvents imply that for every $t>0$, $\displaystyle{\int 
     \frac{1}{t-x}d\nu^{f}_{N}(x)}$ converges to $\displaystyle{\int 
     \frac{1}{t-x}d\nu^{f}(x)}$.\\ By the Stone-Weierstrass theorem, the
    algebra $\mathcal{A}$ of polynomial functions in ${(1-x)^{-1}}$
    defined on
    $\R_{-}$ is dense in the set $C_{\infty}(\R_{-})$ of continuous
    functions on $\R_{-}$ vanishing at $-\infty$. Let $\mathcal{C}$ be
    the vector space spanned by the following set of 
     functions defined on $\R_{-}$, $\{x\mapsto  \frac{1}{t-x},\
    t>0\}$. 
As $\mathcal{A}$ is in the closure of
    $\mathcal{C}$ for the uniform norm,  $\mathcal{C}$ is dense in
    $C_{\infty}(\R_{-})$.
     Thus $$\displaystyle{\int f(x)P_{N,t}^{(n)}(f)(x)dm_{N}^{(n)}(x)=\int_{-\infty}^{0} 
     e^{tx}d\nu^{f}_{N}(x)}$$ converges to $$\displaystyle{\int f(x)P_{t}^{(n)}(f)(x)dm^{(n)}(x)=\int_{-\infty}^{0} e^{tx}d\nu^{f}(x)}$$ for 
     all $t>0$. The polarization identity lets us recover the announced 
    convergence.%
\subsection{Proof  of  theorem \ref{cvnoyaux} }
   First, let us notice that the convergence result stated in 
   proposition \ref{cvmoment} also holds if the invariant measures 
   $m^{(n)}_{N}$ and $m^{(n)}$ of the $n$-point processes are  
   replaced by uniform measures:   
\begin{prop}
\label{convmom}Let $n\in\N^{*}$ and $t>0$.  Let 
    $\lambda_{N,0}$, $\lambda_{N,1}$ and $\lambda$ denote the uniform
    probability  measures on 
    $\frac{1}{N}(\Z/N\Z)$, 
    $\frac{1}{2N}((2\Z+1)/2N\Z)$ and $S^{1}$ respectively. If  $f$ and
    $g$ are  continuous functions on $(S^{1})^n$ then 
    $$\int  
    g(x)P^{(n)}_{N,t}(f)(x)d\lambda_{N,0}^{\otimes n}(x)\mbox{ and }\int g(x)P^{(n)}_{N,t}(f)(x)d\lambda_{N,1}^{\otimes n}(x)$$ converges to 
    $\displaystyle{\int 
    g(x)P^{(n)}_{t}(f)(x)d\lambda^{\otimes n}(x)}$ as $N$ tends to
    $+\infty$.
\end{prop}
\begin{proof}
As $||P^{(n)}_{N,t}(f)||_{\infty}\leq ||f||_{\infty}$ a  density argument
reduces the problem to showing the convergence results for $C^{1}$
functions $f$ and $g$.  
Let $g$ and $f$ be  $C^{1}$ functions on $(S^{1})^{n}$.
First, let us show that $\displaystyle{\int  
    g(x)P^{(n)}_{N,t}(f)(x)d\lambda_{N}^{\otimes n}(x)}$ converges to 
    $\displaystyle{\int g(x)P^{(n)}_{t}(f)(x)d\lambda^{\otimes n}(x)}$ as $N$ tends to
    $+\infty$, where $\lambda^{(n)}_N$ denotes the uniform measure on
    $T^{(n)}_N$. For $\eps\geq 0$, set
$$V_{\eps}=\{x\in (S^{1})^n,\ \exists\ i\neq j,\
      |x_i-x_j|\leq \eps \}$$ and for $\eps>0$, 
consider a continuous function 
$g_{\eps}$  on $(S^{1})^n$ such that $g_{\eps}=0$ on $V_0$,   $g_{\eps}=g$ outside   $V_\eps$
and ${||g_{\eps}||_{\infty}\leq ||g||_{\infty}}$. \\
Outside  $V_0$, the measures  $m_{N}^{(n)}$ and $m^{(n)}$ coincide with 
$\lambda^{(n)}_{N}$ and $\lambda^{\otimes n}$  respectively. It follows from
proposition \ref{cvmoment} that 
$\int g_{\eps}P_{N,t}^{(n)}(f)d\lambda_{N}^{(n)}$ converges
to $\int g_{\eps}P_{t}^{(n)}(f)d\lambda^{\otimes n}$. As
$\lambda^{\otimes n}(\partial
V_{\eps})=0$, $\lambda^{(n)}_{N}(V_{\eps})$ converges to
$\lambda^{\otimes n}(V_{\eps})$. Thus, the upper limit as $N$ tends
to $+\infty$ of 
\begin{multline*}
|\int gP_{N,t}^{(n)}(f)d\lambda^{(n)}_N-\int
gP_{t}^{(n)}(f)d\lambda^{\otimes n}|
\leq \int |g-g_{\eps}||P_{N,t}^{(n)}(f)|d\lambda^{(n)}_N \\
+ |\int g_{\eps}P_{N,t}^{(n)}(f)d\lambda_{N}^{(n)}-\int
g_{\eps}P_{t}^{(n)}(f)d\lambda^{\otimes n}|+
\int|g-g_{\eps}||P_{t}^{(n)}(f)|d\lambda^{\otimes n}
\end{multline*} is bounded by 
$ 4||g||_{\infty}||f||_{\infty}\lambda^{\otimes n}(V_{\eps})$. 
As  $\lambda^{\otimes n}(V_{\eps})$ converges to $\lambda^{\otimes
  n}(V_{0})=0$
as $\eps$ tends to $0$, $|\int gP_{N,t}^{(n)}(f)d\lambda^{(n)}_N-\int
gP_{t}^{(n)}(f)d\lambda^{\otimes n}|$ converges to $0$ as $N$ tends to
$+\infty$.\\
To conclude, let us remark that for every Lipschitz function $\phi$ on
$(S^{1})^n$, the difference between $\int \phi d\lambda^{(n)}_N$ and $\int
\phi d\lambda^{\otimes n}_{N,0}$ or $\int\phi d\lambda^{\otimes n}_{N,1}$ is bounded by 
$$
\frac{1}{2N^{n}}\sum_{u_1=0}^{N-1}\!\!\cdots\!\!\sum_{u_n=0}^{N-1}
|\phi(\frac{u_1}{N},\cdots,\frac{u_n}{N})-\phi(\frac{2u_1+1}{2N},\cdots,\frac{2u_n+1}{2N})|\leq
\frac{n}{4N}||\phi||_{Lip}.$$ As the  Lipschitz coefficient of
$P_{N,t}^{(n)}(f)$  is  bounded by $||f||_{Lip}$, this ends the proof of
proposition \ref{convmom}. 
\end{proof} 
\subsubsection{Proof of  assertion $(i)$}
Let $t$ be a fixed positive real. We shall prove that for every
$j\in\N^*$, if $\mu_N$ for $N\in \N^*$ is a
probability measure defined on $(\frac{1}{N}(\Z/N\Z))^{j}$ such
that $(\mu_N)_N$ converges weakly to a probability measure $\mu$ on
$(S^{1})^{j}$ and if $g$ is a continuous function on $(S^{1})^{j}$
then the random variables $\mu_NK_{N,0,t}^{\otimes j}(g)$
converge in law towards $\mu K_{0,t}^{\otimes j}(g)$ as $N$ tends to
$+\infty$. This convergence result applied  to 
$\mu_N=\nu_N^{\otimes r}$ and $j=\l r$ where  $r\in\N^*$, will show that  
the  measures $E((\nu_NK_{N,0,t}^{\otimes \l})^{\otimes r})$
converge weakly  to $E((\nu K_{0,t}^{\otimes \l})^{\otimes r})$. Thus the
convergence in law of $\nu_N K_{N,0,t}^{\otimes \l}$ will follow. \\
Let  $j\in \N^*$ and let $(\mu_{N})_N$ be a sequence of probability measures
defined on the sets $(\frac{1}{N}(\Z/N\Z))^{j}$ that converges weakly
to a probability measure $\mu$ on $(S^{1})^j$. 
As the $C^1$ functions on $(S^{1})^{j}$ are dense in $C((S^{1})^{j})$, it suffices to prove that for every $g\in C^{1}((S^{1})^{j})$, the
sequence of random variables $(\mu_NK_{N,0,t}^{\otimes j}(g))_N$
converges in law to $\mu K_{0,t}^{\otimes j}(g)$. 
Indeed, let $(g_k)$ be a sequence  of $C^1$ functions  that converges to  $g\in
C((S^{1})^{j})$. For every $u\in\R$, 
\begin{multline*}|E(e^{iu\mu_{N}K_{N,0,t}^{\otimes j}(g)})- E(e^{iu\mu
    K_{0,t}^{\otimes j}(g)})|\leq
|u|E(|\mu_N K_{N,0,t}^{\otimes j}(g_k)-\mu_N K_{N,0,t}^{\otimes
  j}(g)|)\\
+|E(e^{iu\mu_{N}K_{N,0,t}^{\otimes j}(g_k)})-
E(e^{iu\mu K_{0,t}^{\otimes j}(g_k)})|+|u|E(|\mu K_{0,t}^{\otimes
  j}(g_k)-\mu K_{0,t}^{\otimes j}(g)|)\\\leq 2|u|||g-g_k||_{\infty}+
|E(e^{iu\mu_N K_{N,0,t}^{\otimes j}(g_k)})-E(e^{iu\mu K_{0,t}^{\otimes j} (g_k)})|
\end{multline*}
Let $\eps>0$. If we take the upper limit  of the two parts of the previous inequality as $N$ tends to
$+\infty$ with an integer
$k$ satisfying  $||g-g_{k}||_{\infty}\leq \eps$, then  $$\overline{\lim}_{N}|E(e^{iu\mu_{N}K_{N,0,t}^{\otimes j}(g)})- E(e^{iu\mu
    K_{0,t}^{\otimes j}(g)})|\leq 2|u|\eps$$ provided that
$\mu_{N}K_{N,0,t}^{\otimes j}(h)$ converges in law to
$\mu K_{0,t}^{\otimes j}(h)$ for every $h\in C^{1}((S^{1})^{j})$.\\
Let $g$ be a $C^1$ function on $(S^{1})^{j}$. In order to prove that
$X_N=\mu_N K_{N,0,t}^{\otimes j}(g)$ converges in law to $X=\mu K_{0,t}^{\otimes j}(g)$, let us introduce a sequence of probability measures with
Lebesgue density that approaches $\mu$.  
Let $\phi$ be a $C^{\infty}$ density function defined on
$(S^{1})^{j}$. For $k\in \N^*$, set $\phi_{k}: x\mapsto
k^{j}\phi(kx)$. Then any probability measure $\eta$ on $(S^{1})^{j}$ can be approximated by
the probability measures $(\phi_{k}\star \eta)(x) \lambda^{\otimes
  j}(dx)$: more precisely, there is a constant $C_{j,\phi}$ such that
for every Lipschitz function $f$ on $(S^{1})^{j}$ and every
probability measure $\eta$ on $(S^{1})^{j}$
\begin{equation}\label{convol}|\int f(\phi_{k}\star\eta) d\lambda^{\otimes j} -\int f d\eta|\leq
\frac{C_{j,\phi}}{k}||f||_{Lip}.\end{equation}
Indeed, 
\begin{multline*}\Big|\int f(\phi_{k}\star\eta) d\lambda^{\otimes j} -\int f d\eta\Big|\\=\Big|\int\int
f(u)\phi_{k}(u-v)d\lambda^{\otimes j}(u)d\eta(v)-\int
f(v)\Big(\int\phi_{k}(u-v)d\lambda^{\otimes j}(u)\Big)d\eta(v)\Big|\\
\leq
||f||_{Lip}\sum_{i=1}^{j}\int\Big(\int|u_i-v_i|\phi_{k}(u-v)d\lambda^{\otimes
  j}(u)\Big)d\eta(v).
\end{multline*}
To obtain  bound (\ref{convol}), it remains to note that 
$$\int |u_i-v_i|\phi_{k}(u-v)d\lambda^{\otimes j}(u) 
=\frac{1}{k}\int |z_i|\phi(z)d\lambda^{\otimes j}(z).$$
Let us introduce auxiliary variables:  for $k\in \N^*$, set
\begin{alignat*}{2}X^{(k)}&=\int_{(S^{1})^{ j}} K_{0,t}^{\otimes  j}(g)(\phi_k\star\mu)
d\lambda^{\otimes  j},&
X_{N}^{(k)}=\int K_{N,0,t}^{\otimes  j}(g)(\phi_k\star\mu)
d\lambda_{N,0}^{\otimes  j},\\
Y^{(k)}_N&=\int K_{N,0,t}^{\otimes  j}(g)(\phi_k\star\mu_N)
d\lambda_{N,0}^{\otimes  j}.
\end{alignat*}
 For every $\l\in \N^*$, the $\l$-th moments of $X^{(k)}$
and $X_{N}^{(k)}$ are the following: 
\begin{eqnarray*}E((X_{N}^{(k)})^{\l})&=&\int P_{N,t}^{(j\l)}(g^{\otimes \l})(\phi_{k}\star
  \mu)^{\otimes \l}d\lambda_{N,0}^{\otimes \l j},\\ 
  E((X^{(k)})^\l)&=&\int P_{t}^{(j\l)}(g^{\otimes \l})(\phi_{k}\star
  \mu)^{\otimes \l}d\lambda^{\otimes \l j}.
  \end{eqnarray*}
  Thus,  it follows from  proposition \ref{convmom} that for every
$k\in \N^{*}$, the moments of any order of $(X_{N}^{(k)})_N$ converge
to the moments of $X^{(k)}$. As the random variables $X^{(k)}_N$
and $X^{(k)}$ are bounded by $||g||_{\infty}||\phi_k||_{\infty}$,  $(X_{N}^{(k)})_N$ converges in law 
to $X^{(k)}$. In order to deduce   that $(X_N)_N$ converges in law to
$X$, we shall prove the following results: 
\begin{itemize}
\item $E((X^{(k)}-X)^2)$ converges to $0$ as $k$ tends to $+\infty$,
\item $\sup_{N\geq k^{j+2}}E((Y^{(k)}_N-X_N)^2)$ converges to $0$ as $k$ tends
  to $+\infty$,
\item for every $k\in \N^*$, $E((X^{(k)}_N-Y_{N}^{(k)})^2)$ converges to
  $0$ as $N$ tends
  to $+\infty$.
\end{itemize}
These three convergence results are sufficient to deduce that the 
characteristic functions of $X_N$ converge pointwise to the
characteristic function of $X$. Indeed,  for every $u\in \R$ and $k\in \N^*$,
$|E(\exp(iuX_N))-E(\exp(iuX))|$ is bounded by 
\begin{multline*}|E(\exp(iuX_N))-E(\exp(iuX))|\leq 
|u|E(|X_N-Y^{(k)}_N|)+
|u|E(|Y^{(k)}_N-X^{(k)}_N|)\\+
|E(\exp(iuX^{(k)}_N))-E(\exp(iuX^{(k)}))|+|u|E(|X^{(k)}-X|).\end{multline*}
Thus for every $\eps>0$, there exists an integer $k_{\eps}$ such that
for every $N\geq k^{j+2}$, 
\begin{multline*}|E(\exp(iuX_N))-E(\exp(iuX))|\leq
|u|\eps+|u|E(|Y^{(k_{\eps})}_N-X^{(k_{\eps})}_N|)\\
+|E(\exp(iuX^{(k_{\eps})}_N))-E(\exp(iuX^{(k_{\eps})}))|.
\end{multline*} By  taking the upper
limit, as $N$ tends to $+\infty$, of the two terms in this inequality,
we obtain that $|E(\exp(iuX_N))-E(\exp(iuX))|$ converges to $0$. 
\begin{itemize}
\item {\it Study of $E((X^{(k)}-X)^2)$:} first,
  \begin{multline*}E((X^{(k)})^2-X^{(k)}X)=\\\int\Big(\int P_{t}^{(2 j)}(g\otimes
  g)(u,v)(\phi_{k}\star \mu)(u)d\lambda^{\otimes  j}(u)\Big)(\phi_k\star\mu)(v)d\lambda^{\otimes  j}(v)\\-\int\Big(\int P_{t}^{(2j)}(g\otimes
  g)(u,v)(\phi_{k}\star \mu)(u)d\lambda^{\otimes  j}(u)\Big)d\mu(v).\end{multline*} 
As the map $P_{t}^{(2 j)}(g\otimes g)$ is Lipschitz, the map $$v \mapsto \int P_{t}^{(2 j)}(g\otimes
  g)(u,v)(\phi_{k}\star \mu)(u)d\lambda^{\otimes  j}(u)$$ is also Lipschitz
  with  Lipschitz coefficient bounded by $||P_{t}^{(2 j)}(g\otimes
  g)||_{Lip}$. 
Thus, bound (\ref{convol}) gives $|E((X^{(k)})^2-X^{(k)}X)|\leq \frac{C_{ j,\phi}}{k}||P_{t}^{(2 j)}(g\otimes
  g)||_{Lip}$.
Similarly, 
\begin{multline*}|E(X^{(k)}X-X^{2})|=\Big|\int\Big(\int  P_{t}^{(2 j)}(g\otimes
  g)(u,v)d\mu(u)\Big)(\phi_k\star\mu)(v)d\lambda^{\otimes  j}(v)\\-\int\Big(\int  P_{t}^{(2 j)}(g\otimes
  g)(u,v)d\mu(u)\Big)d\mu(v)\Big|\leq \frac{C_{ j,\phi}}{k}||P_{t}^{(2 j)}(g\otimes
  g)||_{Lip}.
\end{multline*} Therefore $E((X^{(k)}-X)^2)$ converges to $0$ as $k$
  tends to $+\infty$.
\item {\it Study of $E((Y^{(k)}_N-X_N)^2)$:} the same splitting as before
  yields that 
\begin{multline*}E((Y^{(k)}_N-X_N)^2)=\int
  (F^{1}_{N}(v)-F^{2}_N(v))(\phi_k\star\mu_N)(v)d\lambda_{N,0}^{\otimes
   j}(v)\\-\int (F^{1}_{N}(v)-F^{2}_N(v))d\mu_{N}(v)
\end{multline*} where 
\begin{eqnarray*}F^{1}_N(v)&=&\int P_{N,t}^{(2 j)}(g\otimes
  g)(u,v)(\phi_{k}\star \mu_N)(u)d\lambda_{N,0}^{\otimes  j}(u)\\
F^{2}_N(v)&=&\int  P_{N,t}^{(2 j)}(g\otimes
  g)(u,v)d\mu_N(u).
\end{eqnarray*} Let us bound  more generally, $$\Delta_{N,k}(f)=\int
  f(\phi_k\star\mu_N)d\lambda_{N,0}^{\otimes
   j}-\int fd\mu_{N}$$ for a Lipschitz function $f$ on
  $(S^{1})^{ j}$. $\Delta_{N,k}(f)$ is the sum of two following terms: 
\begin{eqnarray*}\Delta_{N,k}^{(1)}(f)&=&\int f(\phi_{k}\star
  \mu_N)d\lambda_{N,0}^{\otimes  j}-\int f(\phi_{k}\star
  \mu_N)d\lambda^{\otimes  j}\\ 
\Delta_{N,k}^{(2)}(f)&=&\int f(\phi_{k}\star
  \mu_N)d\lambda^{\otimes  j}-\int fd\mu_N.
\end{eqnarray*} 
It follows from  bound (\ref{convol}) that
  $|\Delta_{N,k}^{(2)}(f)|\leq \frac{C_{j,\phi}}{k}||f||_{Lip}$. 
By developing the  convolution term $\phi_{k}\star \mu_N$,
  the expression of $\Delta^{(1)}_{N,k}(f)$ becomes
  $$\Delta^{(1)}_{N,k}(f)=\!\!\int\!\!\Big(\!\int\!\!
  f(u)\phi_{k}(u-v)d\lambda_{N,0}^{\otimes  j}(u)-\int\!\!
  f(u)\phi_k(u-v)d\lambda^{\otimes  j}(u)\Big)d\mu_N(v).$$ Since for
  every $v\in (S^{1})^{ j}$, $u\mapsto
  f(u)\phi_{k}(u-v)$ is a Lipschitz function with Lipschitz
  coefficient bounded by
  $||f||_{Lip}k^j||\phi||_{\infty}+||f||_{\infty}k^{j+1}||\phi||_{Lip}$, we
  obtain $$|\Delta^{(1)}_{N,k}(f)|\leq
  D_{ j}\frac{k^{j+1}}{N}(\frac{1}{k}||f||_{Lip}||\phi||_{\infty}+||f||_{\infty}||\phi||_{Lip})$$
  where $D_{j}$ is a constant irrespective of $f$
  and $\phi$. In conclusion, there is a constant
  $D_{ j,\phi}$ such that for every
  Lipschitz function $f$ on $(S^{1})^{ j}$ and $k\in \N^*$,
  $$\sup_{N\geq k^{j+2}}|\Delta_{N,k}(f)|\leq
  \frac{D_{ j,\phi}}{k}(||f||_{\infty}+||f||_{Lip}).$$ Since the map
  $P_{N,t}^{(2 j)}(g\otimes g)$ is a Lipschitz function verifying 
$$||P_{N,t}^{(2 j)}(g\otimes g)||_{\infty}\leq ||g\otimes g||_{\infty}
  \mbox{ and } 
  ||P_{N,t}^{(2 j)}(g\otimes g)||_{Lip}\leq ||g\otimes g||_{Lip},$$
  $F^{1}_{N}$ and $F^{2}_{N}$ have the same properties. Thus 
$\sup_{N\geq k^{j+2}}E((Y^{(k)}_N-X_N)^2)$ converges to $0$ as $k$ tends
  to $+\infty$. 
\item {\it Study of $E((Y^{(k)}_{N}-X^{(k)}_N)^2)$:} fix an integer $k\in\N^*$. 
  \begin{multline*}|E((Y^{(k)}_{N})^2-Y^{(k)}_NX^{(k)}_N)|\leq\int\!\!\!\int
  \Big(|P_{N,t}^{(2 j)}(g\otimes g)(u,v)|(\phi_k\star
  \mu_N)(u)\\| (\phi_{k}\star \mu_N)(v)-(\phi_{k}\star
  \mu)(v)|\Big)d\lambda_{N,0}^{\otimes  j}(v)d\lambda_{N,0}^{\otimes
   j}(u)\\\leq k^j||g||^{2}_{\infty}||\phi||_{\infty}\int|(\phi_{k}\star \mu_N)(v)-
(\phi_{k}\star
  \mu)(v)|d\lambda_{N,0}^{\otimes  j}(v)
\end{multline*}
The term $|E((X^{(k)}_{N})^2-Y^{(k)}_NX^{(k)}_N)|$ has the same
bound. \\
The weak convergence of $(\mu_N)_N$ to $\mu$ implies that for every
$u\in (S^{1})^{j}$, $\phi_{k}\star \mu_N(u)$ converges to
$\phi_{k}\star \mu(u)$. As  for every $N$, the map $\phi_{k}\star
\mu_N$ is Lipschitz with 
Lipschitz coefficient  bounded by $k^{j+1}||\phi||_{Lip}$,
${(\phi_{k}\star\mu_N)_N}$ converges uniformly to $\phi_{k}\star
\mu$.
It follows that  ${E((Y^{(k)}_{N}-X^{(k)}_N)^2)}$
 converges to $0$ as $N$ tends to $+\infty$, for every $k\in \N^*$.   
\end{itemize}
  \subsubsection{Proof of  assertion $(ii)$}
 By  $(i)$, $\nu_NK_{N,0,t}^{\otimes  \l}(g)$ converges in law to $\nu
  K_{0,t}^{\otimes  \l}(g)$ as $N$ tends to $+\infty$. These random
  variables are bounded by $||g||_{\infty}$. Thus
  $E(\nu_NK_{N,0,t}^{\otimes  \l}(g))=\int P_{N,0,t}^{(\l)}(g)d\nu_N$
  converges to $E(\nu K_{0,t}^{\otimes  \l}(g))=\int
  P_{0,t}^{(\l)}(g)d\nu$.\\
By splitting $f$ into its positive and its negative part, it suffices
  to deal with  a nonnegative function  $f$. 
   If  $\int fd\nu=0$ then
  $|\int P_{N,t}^{(\l)}(g)fd\nu_N|\leq ||g||_{\infty}\int f d\nu_N$
  tends to $0$. Assume now that $\int fd\nu>0$. Then for $N$ large
  enough, $\int fd\nu_N$ is positive. Let us apply the previous result
  with the sequence of probability measures
  $\tilde{\nu}_N(dx)=\frac{f(x)d\nu_N(x)}{\int fd\nu_N}$:
   $(\tilde{\nu}_N)$ weakly converges to
  $\tilde{\nu}(dx)=\frac{f(x)d\nu(x)}{\int fd\nu}$ whence $\frac{1}{\int f
  d\nu_N}\int P_{N,0,t}^{(\l)}(g)fd\nu_N$ converges to $\frac{1}{\int f
  d\nu}\int P_{0,t}^{(\l)}(g)fd\nu$. As $\int fd\nu_N$ tends to  $\int
  fd\nu$,  $\int P_{N,0,t}^{(\l)}(g)fd\nu_N$ converges to 
$\int P_{0,t}^{(\l)}(g)fd\nu$.
\subsubsection{Proof of  assertion $(iii)$}
Set $\mu_{N,t}=\eta_{N}K_{N,0,t}$ and $\mu_{t}=\eta K_{0,t}$ for 
every $t\in \R$. Let us prove by iteration that for every $r\in\N^*$,
the $r$-dimensional distributions of $\mu_N$ converge weakly to
the $r$-dimensional distributions of $\mu$, that is 
for every  $0=t_1<\ldots <t_r$, the distribution of
$(\mu_{N,t_1},\ldots,\mu_{N,t_r})$ converges to the distribution of
$(\mu_{t_1},\ldots,\mu_{t_r})$.\\
 First, this convergence result holds
for $r=1$. Let $r\in\N^{*}$. Assume that the $r$-dimensional distributions of
$\mu_N$ converge weakly to those of $\mu$. Let ${0=t_1<\ldots<t_{r+1}}$. 
To prove that the law  of $(\mu_{N,t_1},\ldots,\mu_{N,t_{r+1}})$
converges to the distribution of  $(\mu_{t_1},\ldots,\mu_{t_{r+1}})$, it suffices to show that for every 
$k_1,\ldots,k_{r+1}\in\N$, and for every  $g_1\in 
C((S^{1})^{k_1}),\ldots,$ $g_{r+1}\in
C((S^{1})^{k_{r+1}})$,
$E(\mu_{N,t_1}^{\otimes k_1}(g_1)\cdots
\mu_{N,t_{r+1}}^{\otimes k_{r+1}}(g_{r+1}))$
converges to $E(\mu_{t_1}^{\otimes k_1}(g_1)\cdots
\mu_{t_{r+1}}^{\otimes k_{r+1}}(g_{r+1}))$.\\  First, let us note that
for every $N\in\N$, 
$\mu_{N,t_{r+1}}$ is equal to $\mu_{N,t_r}K_{N,t_{r},t_{r+1}}$, the random matrix 
  $K_{N,t_{r},t_{r+1}}$ has the same law as $K_{N,0,t_{r+1}-t_{r}}$
    and is independent of $(\mu_{N,t_1},\ldots,\mu_{N,t_{r}})$.
Thus, 
$$E(\mu_{N,t_1}^{\otimes k_1}(g_1)\cdots
\mu_{N,t_{r+1}}^{\otimes k_{r+1}}(g_{r+1}))
=\int_{(S^{1})^{h}}
P_{N,t_{r+1}-t_{r}}^{(h)}(G)Fd\nu_N$$
where 
 $h=\sum_{i=1}^{r+1}k_i$,  $F$ and $G$ are the maps defined on
$(S^{1})^{h}$ by $$F(x,y)=(g_1\otimes\cdots\otimes
g_{r})(x)\mbox{ and } G(x,y)=g_{r+1}(y)$$ for every $x\in
(S^{1})^{\sum_{i=1}^{r}k_i}$, and $ y\in
(S^{1})^{k_{r+1}}$,
and   
 $\nu_N$ is the probability measure \linebreak[4] ${E(\mu_{N,t_1}^{\otimes k_1}\otimes\cdots\otimes
\mu_{N,t_{r-1}}^{\otimes k_{r-1}}\otimes\mu_{N,t_{r}}^{\otimes
  (k_{r}+k_{r+1})})}$. \\
By the iterative assumption, the distribution of
$(\mu_{N,t_1},\ldots,\mu_{N,t_{r}})$  converges to the distribution of
$(\mu_{t_1},\ldots,\mu_{t_{r}})$. Thus in particular, $(\nu_N)_N$
weakly converges to \linebreak[4] ${\nu=E(\mu_{t_1}^{\otimes k_1}\otimes\cdots\otimes
\mu_{t_{r-1}}^{\otimes k_{r-1}}\otimes\mu_{t_{r}}^{\otimes
  (k_{r}+k_{r+1})})}$.\\ It follows from assertion $(ii)$  that $E(\mu_{N,t_1}^{\otimes k_1}(g_1)\cdots
\mu_{N,t_{r+1}}^{\otimes k_{r+1}}(g_{r+1}))$
converges to \linebreak[4]${E(\mu_{t_1}^{\otimes k_1}(g_1)\cdots
\mu_{t_{r+1}}^{\otimes k_{r+1}}(g_{r+1}))}$, ending the proof of  $(iii)$.

\nocite{stflourPitman,Fukushima}


\begin{thebibliography}{1}

\bibitem{EthierKurtz}
S.N. {\sc Ethier} and T.~G. {\sc Kurtz}.
\newblock {\em Markov processes: Characterization and Convergence}.
\newblock Wiley-Interscience, 1986.

\bibitem{Fukushima}
M.~{\sc Fukushima}, Y.~{\sc Oshima}, and M.~{\sc Takeda}.
\newblock {\em Dirichlet forms and symmetric Markov processes}.
\newblock Walter de Gruyter, 1994.

\bibitem{LeJanRaimond1}
Y.~{\sc Le Jan} and O.~{\sc Raimond}.
\newblock Flows, coalescence and noise.
\newblock math.PR/0203221, To appear in The Annals of Probability.

\bibitem{LeJanRaimond3}
Y.~{\sc Le Jan} and O.~{\sc Raimond}.
\newblock Sticky flows on the circle and their noises.
\newblock To appear in Probab. Theory Relat. Fields.

\bibitem{LeJanRaimond2}
Y.~{\sc Le Jan} and O.~{\sc Raimond}.
\newblock Sticky flows on the circle.
\newblock math.PR/0211387, 2002.

\bibitem{EPitman}
E.~{\sc Pitman}.
\newblock The closest estimates of statistical parameters.
\newblock {\em Proc. Camb. Philos. Soc.}, 33:212--222, 1937.

\bibitem{Pitmanurn}
J.~{\sc Pitman}.
\newblock Some developments of the {Blackwell-MacQueen}.
\newblock In L.S.~Shapley T.S.~Ferguson and J.B. Macqueen, editors, {\em
  Statistics, Probability and Game Theory}, volume~30, pages 245--267. IMS
  Lecture Notes-Monograph, 1996.

\bibitem{stflourPitman}
J.~{\sc Pitman}.
\newblock Combinatorial stochastic processes.
\newblock Saint-Flour lecture notes, July 2002.

\end{thebibliography}
\end{document}